%% file: Manuscript.tex
\documentclass[nonblindrev]{informs3} % current default for manuscript submission

\OneAndAHalfSpacedXI % current default line spacing
\usepackage{afterpage}
\usepackage{caption}
\usepackage{amssymb}
\usepackage{float}
\usepackage{dcolumn}
\usepackage{multirow}
\usepackage{color}
\usepackage{enumitem}
\usepackage{slashbox}
\usepackage{pict2e}
\usepackage{microtype}
\usepackage[noend]{algorithmic}
\usepackage{algorithm}
\usepackage{rotating}
\usepackage{pdflscape}
\usepackage[table]{xcolor}
\usepackage{tabularx, booktabs}
\usepackage{algorithmic}
\usepackage{rotating}
\usepackage{hhline}
\usepackage{url}
\newcolumntype{M}[1]{>{\centering}m{#1}}
\newcolumntype{F}{>{\centering\arraybackslash}m{1.5cm}}
\newcolumntype{C}{>{\centering\arraybackslash}p{2.5em}}
\usepackage{natbib}
 \bibpunct[, ]{(}{)}{,}{a}{}{,}%
 \def\bibfont{\small}%
\TheoremsNumberedThrough     % Preferred (Theorem 1, Lemma 1, Theorem 2)
\EquationsNumberedThrough    % Default: (1), (2), ...
\begin{document}
%%%%%%%%%%%%%%%%

% Outcomment only when entries are known. Otherwise leave as is and 
%   default values will be used.
%\setcounter{page}{1}
%\VOLUME{00}%
%\NO{0}%
%\MONTH{Xxxxx}% (month or a similar seasonal id)
%\YEAR{0000}% e.g., 2005
%\FIRSTPAGE{000}%
%\LASTPAGE{000}%
%\SHORTYEAR{00}% shortened year (two-digit)
%\ISSUE{0000} %
%\LONGFIRSTPAGE{0001} %
%\DOI{10.1287/xxxx.0000.0000}%

% Author's names for the running heads
% Sample depending on the number of authors;
% \RUNAUTHOR{Jones}
% \RUNAUTHOR{Jones and Wilson}
% \RUNAUTHOR{Jones, Miller, and Wilson}
% \RUNAUTHOR{Jones et al.} % for four or more authors
% Enter authors following the given pattern:
\RUNAUTHOR{Matl, Hartl, and Vidal}

% Title or shortened title suitable for running heads. Sample:
% \RUNTITLE{Bundling Information Goods of Decreasing Value}
% Enter the (shortened) title:
\RUNTITLE{Workload Equity in Vehicle Routing: The Impact of Alternative Workload Resources}

% Full title. Sample:
% \TITLE{Bundling Information Goods of Decreasing Value}
% Enter the full title:
\TITLE{Workload Equity in Vehicle Routing: \protect \\ The Impact of Alternative Workload Resources}

% Block of authors and their affiliations starts here:
% NOTE: Authors with same affiliation, if the order of authors allows, 
%   should be entered in ONE field, separated by a comma. 
%   \EMAIL field can be repeated if more than one author
\ARTICLEAUTHORS{%
\AUTHOR{P.~Matl, R.F.~Hartl}
\AFF{University of Vienna, Austria, \EMAIL{piotr.matl@univie.ac.at}, \EMAIL{richard.hartl@univie.ac.at}\URL{}}
\AUTHOR{T.~Vidal}
\AFF{Pontifical Catholic University of Rio de Janeiro, Brazil, \EMAIL{vidalt@inf.puc-rio.br}, \URL{}}
% Enter all authors
} % end of the block

\ABSTRACT{%
In practical vehicle routing problems (VRPs), important non-monetary benefits can be achieved with more balanced operational plans which explicitly consider workload equity. This has motivated practitioners to include a wide variety of balancing criteria in decision support systems, and researchers to examine the properties of these criteria and the trade-offs made when optimizing them. As a result, previous studies have provided a much-needed understanding of how different \textit{equity functions} affect the resulting VRP solutions. However, by focusing exclusively on models which balance tour lengths, a critical aspect has thus far remained unexplored -- namely the impact of the \textit{workload resource} subject to the balancing.

In this work, we generalize previous studies to different workload resources, extend the scope of those analyses to additional aspects of managerial and methodological significance, and reevaluate accordingly previous conclusions and guidelines for formulating a balance criterion. We propose a classification of workload resources and equity functions, and establish which general types of balanced VRP models can lead to unintended optimization outcomes. To explore in greater detail the differences between models satisfying the given guidelines, we conduct an extensive numerical study of 18 alternative balance criteria. We base our observations not only on smaller instances solved to optimality, but also on larger instances solved heuristically, to gauge the extent to which our conclusions hold also for larger VRPs.

Overall, our study counter-balances the focus of previous works and reveals the importance of selecting the right workload resource. Although the marginal cost of equity is low for all examined models, the trade-off structure depends primarily on the workload resource, and less on the equity function. For the same resource, we observe a notable degree of overlap between the solutions found with different equity functions. On the other hand, VRP solutions which are well-balanced with respect to one resource were found to be of poor quality when evaluated in terms of other resources. Finally, we observe that solutions with similar cost and balance tend to exhibit similar solution structure, a property which should be directly exploited by optimization methods.

}%

% Sample
%\KEYWORDS{deterministic inventory theory; infinite linear programming duality; 
%  existence of optimal policies; semi-Markov decision process; cyclic schedule}

% Fill in data. If unknown, outcomment the field
%\HISTORY{}

\maketitle
%%%%%%%%%%%%%%%%%%%%%%%%%%%%%%%%%%%%%%%%%%%%%%%%%%%%%%%%%%%%%%%%%%%%%%

% Samples of sectioning (and labeling) in TRSC
% NOTE: (1) \section and \subsection do NOT end with a period
%       (2) \subsubsection and lower need end punctuation
%       (3) capitalization is as shown (title style).
%
%\section{Introduction.}\label{intro} %%1.
%\subsection{Duality and the Classical EOQ Problem.}\label{class-EOQ} %% 1.1.
%\subsection{Outline.}\label{outline1} %% 1.2.
%\subsubsection{Cyclic Schedules for the General Deterministic SMDP.}
%  \label{cyclic-schedules} %% 1.2.1
%\section{Problem Description.}\label{problemdescription} %% 2.

% Text of your paper here

\input{sections/1-intro}

\input{sections/3-theory}
\input{sections/4-analysis}
\input{sections/5-conclusion}

% Acknowledgments here
%\ACKNOWLEDGMENT{%
% Enter the text of acknowledgments here
%}% Leave this (end of acknowledgment)

% Appendix here
% Options are (1) APPENDIX (with or without general title) or 
%             (2) APPENDICES (if it has more than one unrelated sections)
% Outcomment the appropriate case if necessary
%

%
%   or 
%
% \begin{APPENDICES}
% \section{<Title of Section A>}
% \section{<Title of Section B>}
% etc
% \end{APPENDICES}

% References here (outcomment the appropriate case) 

% CASE 1: BiBTeX used to constantly update the references 
%   (while the paper is being written).
\vspace*{0.3cm}
\bibliographystyle{informs2014trsc} % outcomment this and next line in Case 1
\bibliography{bibliography} % if more than one, comma separated

% CASE 2: BiBTeX used to generate mypaper.bbl (to be further fine tuned)
%\input{mypaper.bbl} % outcomment this line in Case 2

\end{document}

%% file: sections/1-intro.tex
\section{Introduction}
\label{1-intro}

Workload equity and balanced resource utilization are becoming increasingly important issues in real-world logistics systems. This stems from the recognition that logistics is not exclusively cost-driven, and that the marginal cost of more balanced operations can be offset by gains through lower overtime hours, higher employee satisfaction with lower turnover, better customer service, and more flexible use of available capacity. These considerations are relevant in a wide variety of practical vehicle routing problems (VRPs), including service technician routing, vendor-managed inventory systems, waste collection, parcel delivery, public transportation, and volunteer organizations, among others (cf.~the survey in \cite{matl2017a} for a detailed overview and analysis). The heterogeneity of these applications is reflected in the many different forms of ``balance'' objectives which have been proposed.

The implications of choosing one balance objective over another have only recently started to be explored. Abstracting from the particularities of specific applications and considering instead prototypical VRPs, several studies have examined how the structure of optimal VRP solutions changes depending on the chosen balance criterion, in single-objective \citep{campbell2008, huang2012, bertazzi2015} as well as in bi-objective contexts \citep{halvorsen2016, matl2017a}. All of these articles conclude that the choice of a balance criterion has a significant effect on the resulting VRP solutions, and that this choice should therefore not be made arbitrarily.

However, the generality of the above studies is restricted. Specifically, only a particular class of balance criteria has thus far been examined -- namely those balancing tour distances. Yet the articles reporting on practice also describe real-world VRPs in which alternative definitions of workload are more appropriate, e.g.~the number of stops in small package delivery, service time in technician routing and home healthcare services, and load/demand in groceries distribution. By focusing exclusively on the impact of the \textit{equity function} used to quantify the degree of balance, previous studies have not accounted for the choice of the \textit{workload resource} subject to the balancing.

In this article we generalize and extend the previous studies to multiple workload resources, and reassess accordingly previous conclusions and guidelines for formulating a balance objective. We first propose a simple but accurate classification of workload resources as well as equity functions, and identify certain combinations which should be avoided in principle. Since many important aspects for choosing a balance criterion are not amenable to a purely analytical treatment, we also conduct a numerical study of 18 alternative balance objectives based on different combinations of workload resources and equity functions. By considering the bi-objective problem of optimizing both cost and balance, we generate and examine the sets of compromise solutions which would be found in practice with typical constraint-based or weighted-sum approaches.

In the first half of this study we broaden the scope of the research questions examined in \cite{matl2017a} to consider different workload resources, specifically:

\vspace*{0.15cm}
\begin{enumerate}
	\item How is the number of compromise solutions affected by the workload resource and/or the equity function? Do some balanced VRP models yield a wider scope for balancing than others?
	\item What is the marginal cost of balance for different balance objectives? Are there notable differences depending on the workload resource and/or the equity function?
	\item To what extent do different balanced VRP models overlap and generate identical compromise solutions?
\end{enumerate}
\vspace*{0.15cm}

In the second half of our study we further extend earlier analyses by exploring two previously unexamined issues:

\vspace*{0.15cm}
\begin{enumerate}
\setcounter{enumi}{3}
	\item From a managerial perspective: When different balance criteria identify different solutions, does the quality of those solutions deteriorate significantly when evaluated according to a different balancing criterion? In other words, do solutions from different models represent very different preference structures, or is there broad agreement between some models as to what types of VRP solutions are ``well-balanced''?
	\item From a methodological perspective: How strong is the connection between solution similarity in the objective space and similarity in the decision space? To what extent do ``well-balanced'' compromise solutions share a common solution structure that can be exploited by search algorithms?
\end{enumerate}
\vspace*{0.15cm}

Finally, we conduct our analysis separately on smaller instances solved to optimality, and then on larger instances solved with an adaptation of a state-of-the-art VRP heuristic \citep{matl2017b}. By comparing and contrasting the computational results on these two datasets, we investigate and quantify the degree to which previous observations and conclusions still hold as instance size grows, which has thus far not been considered by any of the previous studies.

%% file: sections/3-theory.tex
\section{Equity Objectives and their Properties}
\label{3-theory}

When dealing with balancing concerns, it is important to realize that the equitable distribution of resources is in itself already a multi-objective problem. Distributing a resource between $n$ agents (e.g.~drivers, territories, machines) corresponds to an $n$-objective problem in which each objective $n_i$ represents the allocation to a different agent $i$. From this perspective, it is clear that any equitable allocation must be a Pareto-optimal solution to this $n$-objective problem \citep{kostreva2004}. In practice however, the prevalence of non-dominated or incomparable solutions increases rapidly even for small $n$, and so the classical Pareto-dominance relation is of limited help in ranking different allocations. Hence the need for inequality measures which reduce the dimensionality of equity and provide a stricter ordering of alternative resource distributions.

Inequality measures are designed to capture the equity of a resource allocation with a single number. Formally, an inequality measure $I(\mathbf{x})$ transforms a vector $\mathbf{x}$ of outcomes or workloads $\{x_1, x_2, ..., x_n\}$ into a scalar. In this way, the resulting equity criterion may be expressed in the form of a more standard objective or constraint. In the following, we will distinguish two essential components of an inequality measure: the \textit{resource} (or metric) that determines the workloads in $\mathbf{x}$, and the \textit{function} that transforms this vector into a scalar.

\subsection{Workload Resources}
\label{3-1}

Despite the heterogeneous nature of the balanced VRPs presented in the literature, their models fall into three categories of workload resources: those balancing the \textit{distance/duration} of the tours, those based on the \textit{demand/load}, and those counting the number of \textit{stops/customers}. These alternatives can be generalized into more abstract groups: those in which the workload occurs on the \textit{edges} of the network (e.g.~distance), at the \textit{nodes} (e.g.~demand), or on both (e.g.~tour duration, a function of distance and service time). This reveals that in a typical VRP, the stops/customers resource is a special case of demand/load.

However, the distinction between edge-based and node-based resources can be misleading, since the meaning of edges and nodes can be arbitrary -- consider for instance VRPs and arc routing problems \citep{pearn1987} -- and other problems might be more intuitively formulated without the explicit use of network or graph structures. Rather, the essential distinction to be made is whether the total workload to be distributed is constant for all feasible solutions to the problem, or whether it can vary from one solution to another.

\vspace*{0.3cm}
\begin{definition}
\label{def1}
A workload resource is \textbf{constant-sum} if $\sum_{i=1}^{n} x_i$ is identical for all feasible workload allocations $\mathbf{x} \in \mathbf{X}$, and \textbf{variable-sum} otherwise.
\end{definition}
\vspace*{0.3cm}

Recall that an equitable allocation should be a Pareto-optimal solution to the $n$-objective allocation problem. It is easy to verify that any constant-sum resource satisfies this property -- any change in one outcome corresponds to a change of equal magnitude in another outcome, leading to a non-dominated alternative allocation. In contrast, the more general category of variable-sum resources does not guarantee the satisfaction of the Pareto-optimality principle. Either way, a suitable equity function still needs to be selected in order to introduce a stricter ordering among all the potential allocations. For variable-sum resources, the equity function should additionally prevent the acceptance of dominated alternatives.

%For both types of resources, a suitable equity function still needs to be selected in order to introduce a stricter ordering among all the potential allocations. For constant-sum resources, it is critical since every allocation is Pareto-optimal and thus this relation is of little help, and for variable-sum resources, the equity function should additionally prevent the acceptance of dominated alternatives.

\subsection{Equity Functions}
\label{3-2}

There exists a large body of research on equity functions in the economics and sociology literature, predominantly in the fields of income inequality and wealth distribution (see e.g.~\citealt{sen1973, allison1978}). It is important to note that from a theoretical perspective, there is no consensus on a ``perfect" equity function -- even by approaching this issue from an axiom-based direction, there is disagreement about which axioms should or should not be considered, and this tends to vary depending on the intended area of application. This is not surprising -- and perhaps unavoidable -- given the inherently multi-objective nature of equity and fairness.

Under these circumstances, it is not our intention here to identify an ideal equity function, nor is it to review all potential issues surrounding any specific function. Instead, we will focus our attention on two properties that are, in our view, particularly critical in practice and which \textit{interact directly} with the type of workload resource. A review of other axioms can be found in \cite{matl2017a}.

\begin{axiom}
\textbf{Monotonicity}: Let $\mathbf{x'}$ be formed as follows: $x_i = x_i + \delta_i$ for all $i$ in $\mathbf{x}$. If $\delta_i \geq 0$ for all $i$ with at least one strict inequality, then $I(\mathbf{x'}) \geq I(\mathbf{x})$.
\end{axiom}

The monotonicity axiom is closely related to the Pareto-optimality principle discussed in the previous section. It states that if one or more outcomes are strictly worsened while leaving all others unchanged, then inequality should increase (strong monotonicity) or at least not be reduced (weak monotonicity). In effect, this penalizes Pareto-dominated allocations. As a result, monotonic equity functions minimize inequality in a way that is consistent also with the minimization of individual outcomes \citep{ogryczak2014}. 

\begin{axiom}
\textbf{Transfer Principle}: Let $\mathbf{x'}$ be formed as follows: $x'_i = x_i - \delta$, $x'_j = x_j + \delta$, $x'_k = x_k$, for all $k \notin \{i,j\}$. If $\delta > 0$ and $x_j \geq x_i$, then $I(\mathbf{x'}) \geq I(\mathbf{x})$.
\end{axiom}

The Pigou-Dalton (PD) principle of transfers is one of the most commonly accepted axioms for equity functions. Originating in the context of income distribution, the transfer principle is expressed in terms of maximization objectives, but minimization objectives such as workload can be transformed by subtracting from a sufficiently large constant (i.e.~$x'_i = M - x_i$).

Informally, the transfer principle states that taking from the poorer and giving to the richer should lead to a numerical increase in inequality (strong version), or at least not reduce it (weak version). More formally, a transfer refers to shifting $\delta$ units from an individual $i$ to another individual $j$, such that $x'_i = x_i - \delta$ and $x'_j = x_j + \delta$. A transfer is said to be \textit{regressive} (favoring the party who is already better off) if $x_i \leq x_j$, and \textit{progressive} if $0 \leq \delta \leq x_i - x_j$. If an allocation $\mathbf{x'}$ can be reached by a finite series of only regressive transfers from an allocation $\mathbf{x}$, then $I(\mathbf{x'}) > I(\mathbf{x})$ (strong version), or at least $I(\mathbf{x'}) \geq I(\mathbf{x})$ (weak version). Note that the transfer principle is ambiguous between allocations connected by a combination of regressive \textit{and} progressive transfers, and it applies only if the allocations do not differ in the number and the sum of their outcomes.

\subsection{Implications}
\label{3-3}

The formulation of an equity or balance objective (an inequality measure) implies the selection of both a workload resource as well as an equity function -- the question thus arises about the potential interactions between these two decisions. In that regard, we can make the following observations:

\vspace*{0.3cm}
\begin{enumerate}
	\item For the \textit{special case} of constant-sum resources, the Pareto-optimality of every allocation is guaranteed regardless of the selected equity function. The monotonicity property is therefore of no consequence for constant-sum resources, and so even non-monotonic functions can be used freely. In contrast, the transfer principle is much more relevant, because the constant-sum character of the resource guarantees that all allocations are connected by mean-preserving transfers and thus comparable.
	%\newpage
	\item For the \textit{general case} of variable-sum resources, the Pareto-optimality of every allocation is not guaranteed, and so care must be taken to avoid an active search for and acceptance of dominated solutions to the $n$-objective allocation problem. This makes the monotonicity property critical for the choice of the equity function -- a non-monotonic equity function may imply a preference even for those allocations in which all outcomes are worse than those of another feasible allocation. On the other hand, the transfer principle is of only limited importance: it will apply only in the comparatively rare case when two competing allocations have the same sum of outcomes.
\end{enumerate}
\vspace*{0.3cm}

These guidelines are summarized in Table \ref{guidelines}.

\input{tables/guidelines}

Balance objectives which do not satisfy these criteria can lead to pathological optimization outcomes. This is especially the case for the combination of a variable-sum resource and a non-monotonic function -- which is unfortunately the most common combination, as it places no particular requirements on the function and variable-sum resources represent the general case. In routing problems such as the VRP with route balancing (VRPRB), combining tour lengths with non-monotonic functions can lead to optimal solutions with non-TSP-optimal tours \citep{jozefowiez2002} or even TSP-optimal solutions whose tours are all longer than in other feasible solutions \citep{matl2017a}. But similar observations can be made also for other problem classes. For example in facility location problems, using non-monotonic functions to balance customers' distance to their nearest facility leads to solutions in which all facilities are located far away, providing an equal \textit{lack} of service \citep{marsh1994, eiselt1995, ogryczak2000}. In a scheduling environment, balancing the makespans of individual machines according to a non-monotonic function could lead to plans with artificially long idle times but perfectly equal makespans. Outcomes such as these clearly undermine the credibility, usefulness, and acceptance of the models, methods, and solutions proposed by researchers. Balance objectives should therefore be chosen carefully.

Although the presented guidelines are general, models which satisfy them can still differ greatly depending on the particular workload resource and equity function used. Many important factors for choosing a balance objective -- such as the number of potential trade-off solutions, the marginal cost of improving balance, or the effect on solution structure -- defy meaningful description exclusively through analytical means. We have therefore conducted an extensive numerical study to examine how the choice of the workload resource, the equity function, and their combinations, affect the obtained VRP solutions. In the next section we describe in detail our experiments and report the most relevant observations made from the study.

%% file: tables/guidelines.tex
%\begin{table}[]
%\centering
%\begin{tabular}{c|c|c}
%\textbf{Workload Resource} & \textbf{Monotonicity} & \textbf{Transfer Principle} \\ \hline
%constant-sum & $\medcircle \medcircle \medcircle$ & $\medblackcircle \medblackcircle \medcircle$ \\
%variable-sum & $\medblackcircle \medblackcircle \medblackcircle$ & $\medblackcircle \medcircle \medcircle$
%\end{tabular}
%\caption{Importance of an equity function's satisfaction of the monotonicity axiom and the PD transfer principle, according to the type of workload resource. Importance ranges from irrelevant ($\medcircle \medcircle \medcircle$), to marginal, useful, and critical ($\medblackcircle \medblackcircle \medblackcircle$).}
%\label{guidelines}
%\end{table}

\begin{table}[]
\centering
\begin{tabular}{c|c|c}
\textbf{Workload Resource} & \textbf{Monotonicity} & \textbf{Transfer Principle} \\ \hline
constant-sum & irrelevant & useful \\
variable-sum & critical & marginal
\end{tabular}
\caption{Importance of an equity function's satisfaction of the monotonicity axiom and the Pigou-Dalton transfer principle, according to the type of workload resource.}
\label{guidelines}
\end{table}

%% file: sections/4-analysis.tex
\section{Numerical Study}
\label{4-results}

%Short intro text.

\input{sections/4-1-study-design}
\input{sections/4-2-cardinality}

\input{sections/4-3-cost}

\input{sections/4-4-overlap}
\input{sections/4-5-agreement}
\input{sections/4-6-similarity}

%% file: sections/4-1-study-design.tex
%\subsection{Experimental Design}
%\label{4-1}

We examine in turn the three workload resources identified in Section \ref{3-1}: the distance/duration, the demand/load, and the number of stops/customers served per tour. For each of these resources, we consider six commonly used equity functions: minimization of the maximum workload, lexicographic minimization, the range, the mean absolute deviation, the standard deviation, and the Gini coefficient. This yields a total of 18 alternative objectives for the balancing criterion.

\vspace*{0.3cm}
\paragraph{\textbf{Bi-objective Models.}}
Equity and balancing are typically not the primary optimization objective in VRPs. However, the inclusion of an equity or balance constraint in standard single-objective models implies that the resulting solution will be (in the best case) one of the Pareto-optimal compromise solutions to the bi-objective extension of the original problem. The proper choice of a balancing criterion -- and especially the specific constraint value -- is therefore just as impactful as the primary optimization objective. Hence, we consider the bi-objective problem of minimizing both cost and imbalance, and analyze the resulting Pareto-optimal or non-dominated solution sets obtained with each of the 18 alternative models. This allows us to make broader observations about the models and to make statements about the compromise solutions in general, rather than only about the objective function extremes or some other arbitrarily chosen solutions.

%\subsubsection{Instances}

\vspace*{0.3cm}
\paragraph{\textbf{Instances.}}
Our numerical study is based on two groups of instances. The first group consists of small instances of the Capacitated VRP (CVRP) ($n=20$ customers, $K=5$ vehicles). A total of 40 such instances were generated based on the CVRP benchmark of \cite{uchoa2017}: by taking the first, second, third, etc.~set of 21 customers and labeling the first customer of each set as the depot, we created 20 instances with various customer-depot configurations from the data of instance X62 (demands $q_i$ from $[50,100]$, i.e.~small variance), and another 20 small instances with the data of instance X64 (demand values in $[1,100]$, i.e.~large variance). The vehicle capacity was set to $Q = \lceil \frac{1}{K-1} \sum_{i=1}^{n} q_i \rceil - 1$ so that all vehicles are required while allowing some excess capacity for balancing. We solve this set of instances to optimality by means of enumeration.

The second group of instances consists of the benchmark sets A and B introduced by \cite{augerat1995} and the first 20 instances with non-unitary demands from the X benchmark of \cite{uchoa2017} (unitary instances were disregarded as they would produce identical solutions for the load and stops resources). This yields a total of 70 diverse instances with between 30 and 200 customers, 5 to 50 vehicles, random, clustered, and mixed customer locations, as well as various demand distributions. %All of these instances, as well as their corresponding cost-optimal solutions, are available on the \citeauthor{cvrplib} repository.

We solve the larger instances using the Iterated Hybrid Genetic Search (IHGS) described in \cite{matl2017b}. IHGS significantly outperforms the current state-of-the-art on the VRPRB benchmark. In addition, we have validated its performance on all 18 alternative balance objectives using the 40 smaller instances solved to optimality: IHGS was able to identify over 95\% of the approx.~18,000 Pareto-optimal solutions to these instances. In order to identify solutions of as high a quality as possible also for larger instances, we deactivated the neighborhood granularity parameter, quadrupled the number of iterations per sub-problem, and solved each instance 10 times. We then generated the best known reference set for each model and instance by forming the non-dominated union of all solutions identified in all test runs on that instance and the cost-optimal solutions from the literature (available on the \citeauthor{cvrplib} repository).

The instances used in this study and their optimal/best-known solution sets are included with the online version of this article, and also available upon request from the authors.

\vspace*{0.3cm}
\paragraph{\textbf{Tour Structure.}}
As mentioned in Section \ref{3-3}, combining a variable-sum workload metric with a non-monotonic equity function generally leads to Pareto-optimal VRP solutions with non-TSP-optimal tours \citep{jozefowiez2002}, or TSP-optimal tours whose lengths are all longer than those of other Pareto-optimal solutions, i.e.~total cost is higher, all workloads are worse, yet the solution is deemed Pareto-efficient \citep{matl2017b}. Clearly this contradicts the intent of balanced VRP models. In our study, this affects those models combining the distance resource with a non-monotonic equity function. In those cases, we therefore enforce a TSP-optimality constraint on the instances solved to optimality, and \textsc{move}, \textsc{swap}, and \textsc{2-opt} neighborhood optimality constraints on the instances solved heuristically. In addition, the generated Pareto sets contain only solutions which are workload consistent, i.e. if total cost is higher than in another solution, at least one workload must be strictly lower, otherwise the solution is discarded.

%% file: sections/4-2-cardinality.tex
\subsection{Number of Compromise Solutions}
\label{4-2}

We start by first considering the cardinality of the Pareto sets generated by different balanced VRP models. Pareto set cardinality has an impact on the usefulness of bi-objective models for decision-makers. Generally, a larger number of compromise solutions increases the likelihood that a solution which reflects the decision-maker's true preferences actually exists, since (ceteris paribus) a larger part of the objective space will contain corresponding solutions. This is especially critical when balance considerations are handled with constraints -- the more solutions there are, the higher the likelihood that one exists close to the specified constraint value. Of course, Pareto set cardinality is also a factor in the appropriateness and efficiency of corresponding optimization methods.

 %The data is reported separately for the smaller instances solved to optimality, and for the larger ones solved heuristically.

\begin{figure} [b]
	\centering
		\includegraphics[width=0.95\textwidth]{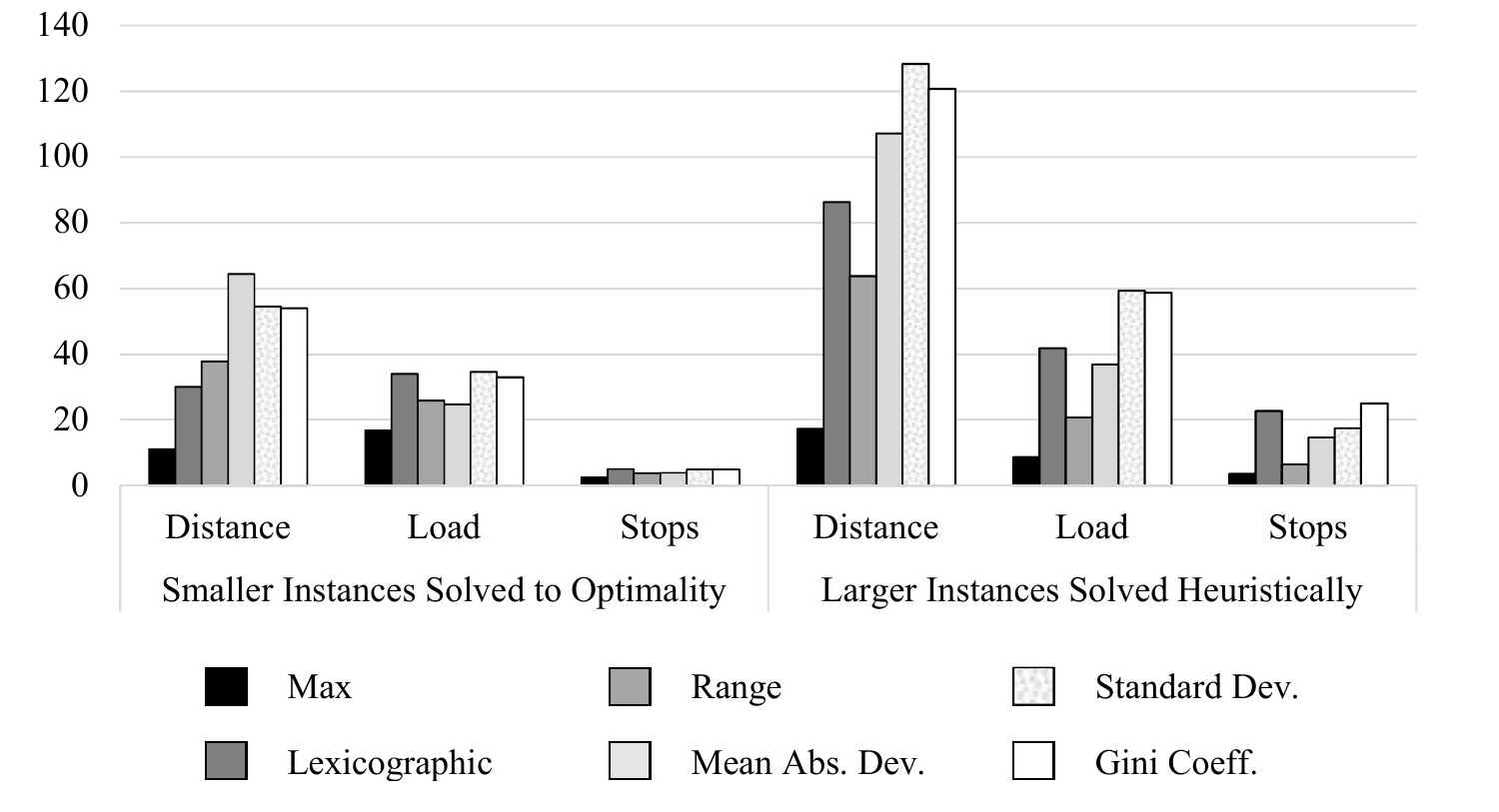}
	\caption{Average solution set cardinality for various balanced VRP models.}
	\label{fig:cardinality}
\end{figure}

Figure \ref{fig:cardinality} visualizes the average Pareto set cardinality according to the chosen workload resource and equity function. One of the first observations we can make is that regardless of the instance size, models balancing distance tend to generate the most trade-off solutions, while those balancing stops have the fewest, with load-based models lying in between. Of course, the larger instances usually lead to larger Pareto sets, but given the size of the smaller instances, the cardinality of their Pareto sets is surprisingly large for the distance- and load-based models. %On the other hand, the models balancing the number of stops provide only quite limited flexibility for balancing, and this flexibility grows rather slowly with instance size.%, especially compared to models balancing the other two resources.

These differences with respect to the workload resource most likely stem from the combinatorial nature of the possible allocations of these resources: the number of different tour lengths is much higher than the number of unique tour loads, which in turn is much higher than the number of tour stops. As a result, fewer unique allocations exist for load and stops, which directly leads to fewer unique balance values and limits the potential for improvement.

When it comes to the equity functions used to measure the balance of different allocations, we can see that the simpler functions (max, range) usually generate smaller Pareto sets, whereas the more complex functions (e.g.~lexicographic, standard deviation) generate larger sets. These differences become more pronounced on the larger instances. This is particularly visible with the outstandingly small Pareto sets generated by the max function, regardless of the resource.%even with the distance-based models which otherwise generate the largest sets.

The explanation for these differences lies most likely in the varying numerical ranges and precision of these equity functions: the number of unique function values places an upper bound on the number of Pareto-optimal trade-off solutions. For example, every workload allocation has a unique lexicographic rank, but millions of allocations may share the same range of workloads. Furthermore, the number of unique values returned by functions like the standard deviation and the Gini coefficient is affected by the number of vehicles, unlike the values taken by the max and range functions. %This interaction becomes especially apparent with the max function, where even using the distance-based models with the largest number of unique workload allocations, the max function identifies fewer than 20 trade-off solutions on average even on the larger instances. 

Overall, we find that richer sets of compromise solutions can be found with the distance and load resources and the more sophisticated equity functions. If a combination of these resources and functions corresponds to a meaningful definition of workload in a given application case, then such models can potentially identify better opportunities to improve balance.

%% file: sections/4-3-cost.tex
\subsection{Cost of Equity}
\label{4-3}

In the preceding section we observed that certain combinations of workload resources and equity functions can lead to larger numbers of compromise solutions, offering a greater degree of opportunity for improving balance. However, the question of whether these opportunities are worthwhile depends critically on the relative cost of improving balance. After all, cost-oriented factors are likely to remain the primary optimization objective in most real-life logistics systems, and better balanced operations are usually of interest if their marginal costs are low. In this section we therefore analyze more closely the typical distribution in the objective space of the trade-off solutions identified by different balanced VRP models.

Figure \ref{fig:sets1} plots for each combination of workload resource and equity function all the Pareto-optimal solutions identified for each of the smaller instances. In this way, both the general trend as well as outliers are visible and comparable on the same plot. Since equity functions quantify the degree of \textit{imbalance}, the balance objective is modeled as a minimization problem. The balance objective function values (y-axis) have been normalized to the range $[0,1]$ in order to make balance improvements (i.e.~imbalance reductions) comparable across different functions and different instances. Similarly, cost values (x-axis) are reported as increases relative to the corresponding cost-optimum (known for all instances in our study).

Figure \ref{fig:sets2} plots the same data as Figure \ref{fig:sets1}, but for the larger instances solved heuristically. Although the optimally balanced solutions are not known, in many cases our heuristic was able to identify solutions with the best possible balance value (e.g.~0), so they are at least weakly Pareto-efficient. 

%\newgeometry{margin=1cm}
\begin{landscape}
\centering\vspace*{\fill}
\thispagestyle{empty} %% Remove header and footer.
\begin{figure} [!htb]
	\centering
		\includegraphics[width=1.0\linewidth]{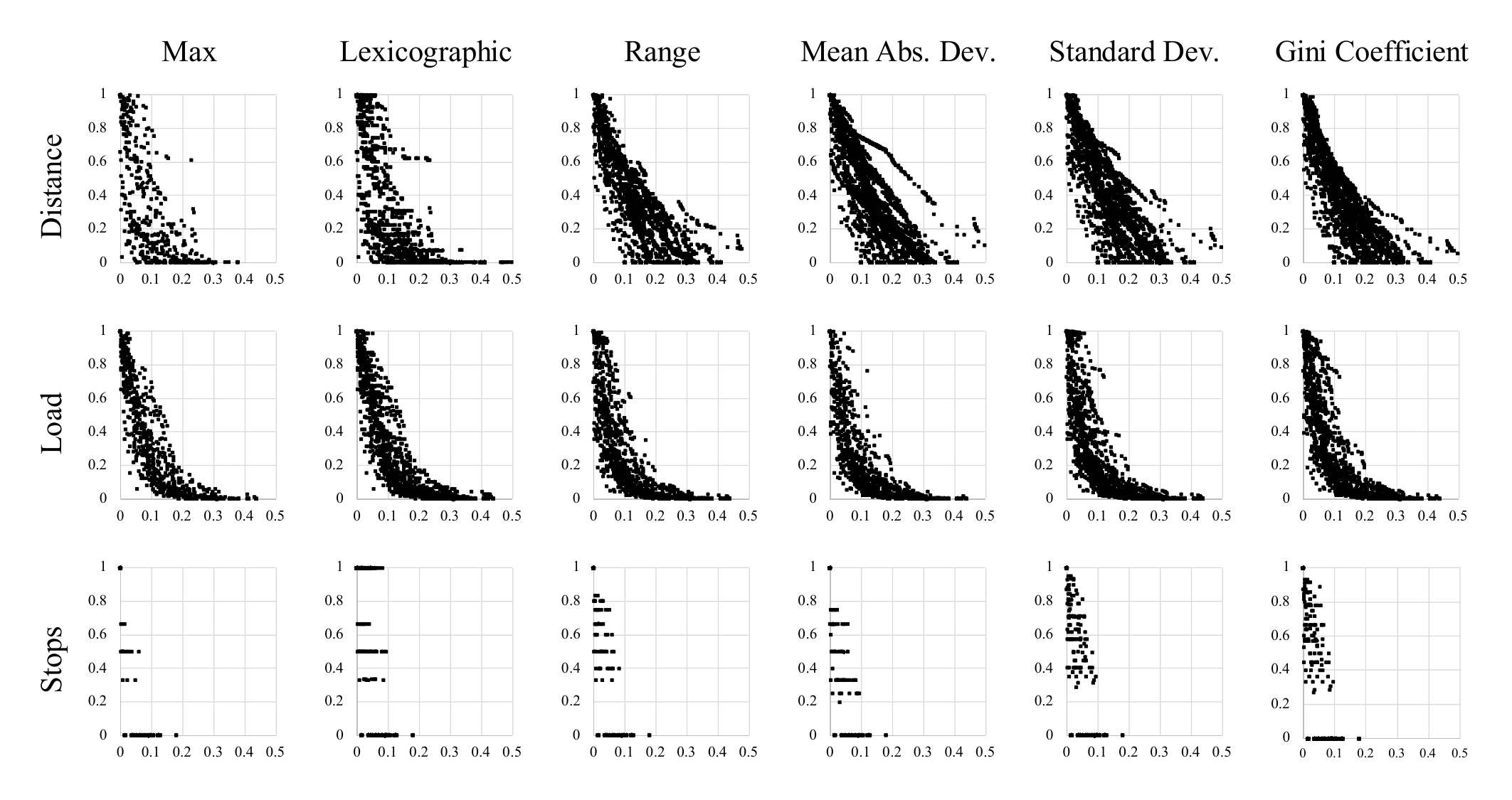}
	\caption{Smaller instances solved to optimality -- distribution of Pareto-optimal solutions according to workload resource and equity function. Cost is plotted on the x-axis as the relative increase above the optimal cost. Balance (modeled as the minimization of imbalance) is plotted along the y-axis and normalized to the interval $[0,1]$ in order to make values comparable across different functions and different instances. Each plot contains all solutions of the corresponding model over all instances in order to visualize the overall trend while also showing outliers.}
	\label{fig:sets1}
\end{figure}
\vfill
\end{landscape}
%\restoregeometry

%\newgeometry{margin=1cm}
\begin{landscape}
\centering\vspace*{\fill}
\thispagestyle{empty} %% Remove header and footer.
\begin{figure} [!htb]
	\centering
		\includegraphics[width=1.0\linewidth]{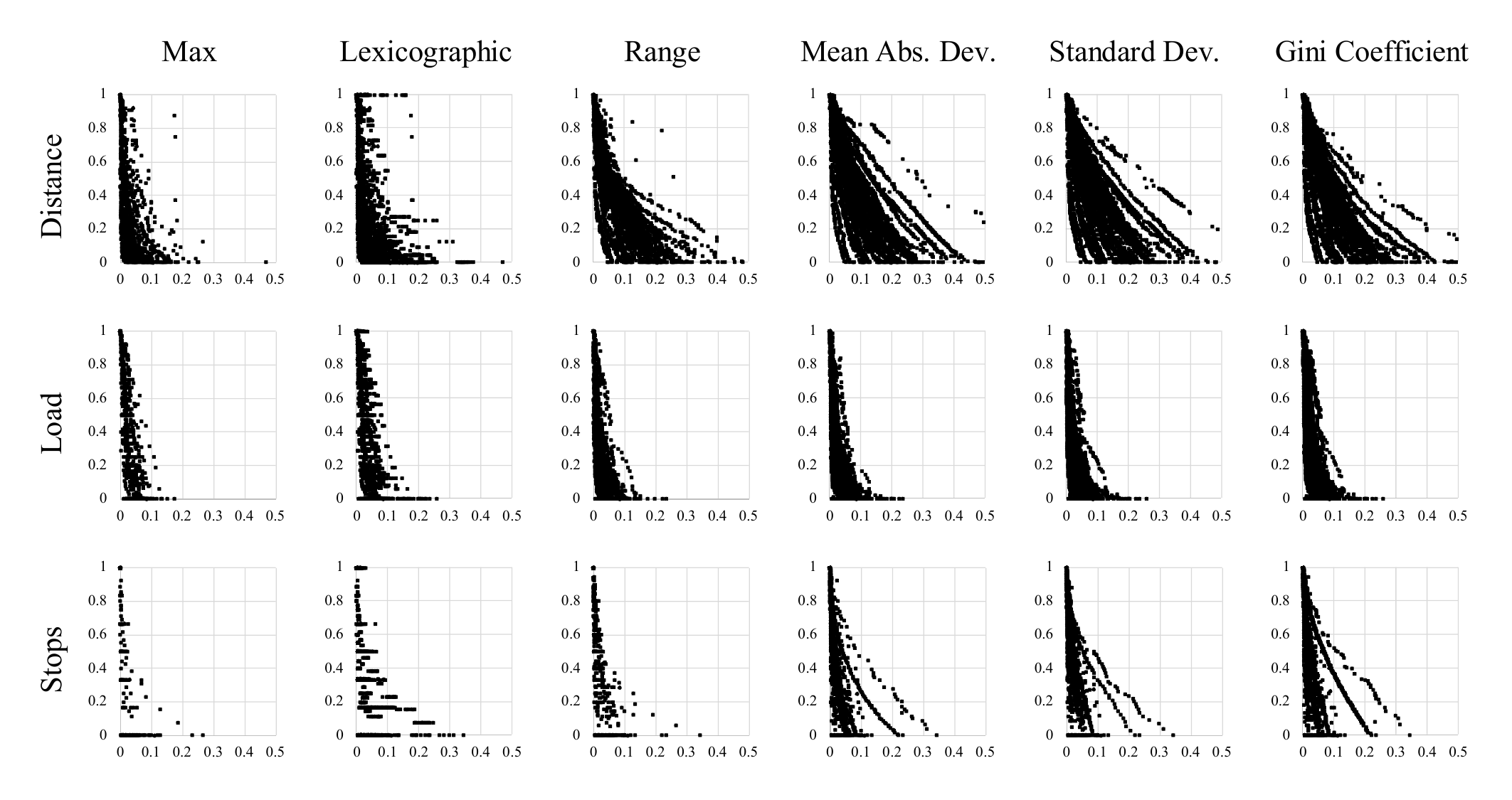}
	\caption{Larger instances solved heuristically -- distribution of best known non-dominated solutions according to workload resource and equity function. Cost is plotted on the x-axis as the relative increase above the optimal cost. Balance (modeled as the minimization of imbalance) is plotted along the y-axis and normalized to the interval $[0,1]$ in order to make values comparable across different functions and different instances. Each plot contains all solutions of the corresponding model over all instances in order to visualize the overall trend while also showing outliers.}
	\label{fig:sets2}
\end{figure}
\vfill
\end{landscape}
%\restoregeometry
%\onehalfspacing

Considering first the smaller instances solved to optimality, we can see that the cost vs.~balance trade-off structure depends primarily on the chosen workload resource, and not on the equity function. For any of the three workload resources, the trade-off structure is remarkably similar for all six equity functions. On the other hand, the differences between models using different workload resources are more pronounced. First, we observe that the solutions of the distance-based models have a more variable distribution in the objective space, meaning that the actual trade-off structure tends to vary more from instance to instance. For the load-based models, we can see that all the Pareto sets combined fall within a narrower region of the objective space, and for stops-based models, the concentration is narrower still. A second observation is that cost and equity appear to be more conflicting objectives when balancing distance-based metrics -- solutions balancing distance tend to be further from the ideal point at the origin than solutions balancing load or stops.

We can identify the same trends also on the larger instances solved heuristically. Here too, the trade-off structure is determined more by the workload resource than by the equity function, and the variability of that structure is also higher for distance-based models than for those balancing load or stops. However, one important difference on the larger instances is that the cost and balance objectives become less conflicting, with most solutions being relatively low-cost, as well as closer to the origin. This suggests that the greater degree of flexibility on larger instances allows for more cost-efficient workload balancing, and overall this appears to hold for all the considered models.

Despite some of these differences, all the considered models have one important aspect in common: the marginal cost of balance is low. For nearly all the models and instances examined in our study, balance can be improved dramatically at little additional cost -- especially on the larger instances, near-optimal balance can be achieved within 10\% or even 5\% of the corresponding cost optimum. Although the actual trade-off structure depends to some extent on the workload resource, all of them allow significant balance improvement at comparatively low cost, and this potential grows with instance size. From a managerial perspective, this allows to more easily take advantage of the various non-monetary benefits of balanced operations discussed in the literature. From a methodological perspective, we can expect that typical weighted-sum approaches need not place overly large weights on balance to find noticeable improvements, and constraint-based methods can use relatively ambitious balance targets without impacting cost severely.

%% file: sections/4-4-overlap.tex
\subsection{Overlap Between Different Models}
\label{4-4}

In the previous section we observed that for a given workload resource, the choice of the equity function appears to have limited impact on the overall trade-off structure. The question thus arises to what extent these different functions actually produce different solutions. %Since the various equity functions are used to quantify the same concept, it seems intuitive that there might be at least some degree of consensus between them as to which solutions are ``well-balanced'', and which are not.

\input{tables/sample-sets}

To explore this further, we first present in Table \ref{sample-sets} the Pareto-optimal solution sets for a typical instance balancing tour loads. For each solution we list its total cost, its corresponding balance according to each of the six equity functions, and its workload allocations. A dash indicates that the solution is not Pareto-optimal for the respective equity function. The lexicographic objective has been replaced with a ranking since the respective workload vector is given on the right.

Considering this sample instance, we can see that there are relatively few solutions which are Pareto-optimal for \textit{every} examined function. Excluding the cost optimum, only seven other solutions fall into this category, out of over 60 total. At the same time, there exist solutions which are unique to only a single function: e.g.~solutions 14, 39, and 50 are identified only by the lexicographic objective, over 10 solutions are unique to the range function, solutions 24 and 45 are found only with the mean absolute deviation, solution 53 only with standard deviation, and solution 47 only with the Gini coefficient. Furthermore, even for the same total cost, alternative workload allocations can be preferred by different functions, e.g.~there is no obvious agreement when it comes to solutions 3 and 4, 10 and 11, 23 and 24, 42 and 43, or 46 and 47. We note also the special case of solutions 3 and 4, where the max function can be indifferent to a lexicographically better allocation with the same maximum workload.

Despite the limited occurrence of solutions identified by \textit{every} examined function, there are generally more solutions shared between different objectives than unique to any single one. For example, the six individual Pareto sets in Table \ref{sample-sets} have a combined cardinality of 200, but after removing duplicates only 66 of these solutions are unique. Except for the max function, the other Pareto sets have between 31 and 41 solutions, which implies a considerable degree of overlap if only 66 different solutions are identified overall.

In order to make more general observations, we have conducted the above analysis for each of the workload resources and all instances. Specifically, we consider for each workload resource the combined set of all Pareto-optimal solutions found with all six equity functions over all instances. For each equity function $F$, this set of solutions can be divided into four categories:

\vspace*{0.25cm}
\begin{enumerate} [label=(\Alph*)]
\itemsep0em 
	\item solutions found by $F$ and all other examined functions,
	\item solutions common to $F$ and a strict subset of the other functions,
	\item solutions identified only with $F$, and
	\item solutions not found by $F$.
\end{enumerate}
\vspace*{0.25cm}

Figure \ref{overlap} presents this data for the sets of smaller and larger instances.

We can see that in general, full overlap between all six functions is uncommon even on the smaller instances, and virtually disappears on the larger ones. The share of solutions unique to any particular function is likewise relatively small, but this tends to increase with instance size as larger Pareto sets and a larger solution space widen the scope for unique workload allocations. As expected, a larger share of solutions lies between these two extremes, being found by some, but not all functions.

\begin{landscape}
\centering
\vspace*{\fill}
\thispagestyle{empty} %% Remove header and footer.
\begin{figure} [!htb]
	\centering
		\includegraphics[width=1.5\textwidth]{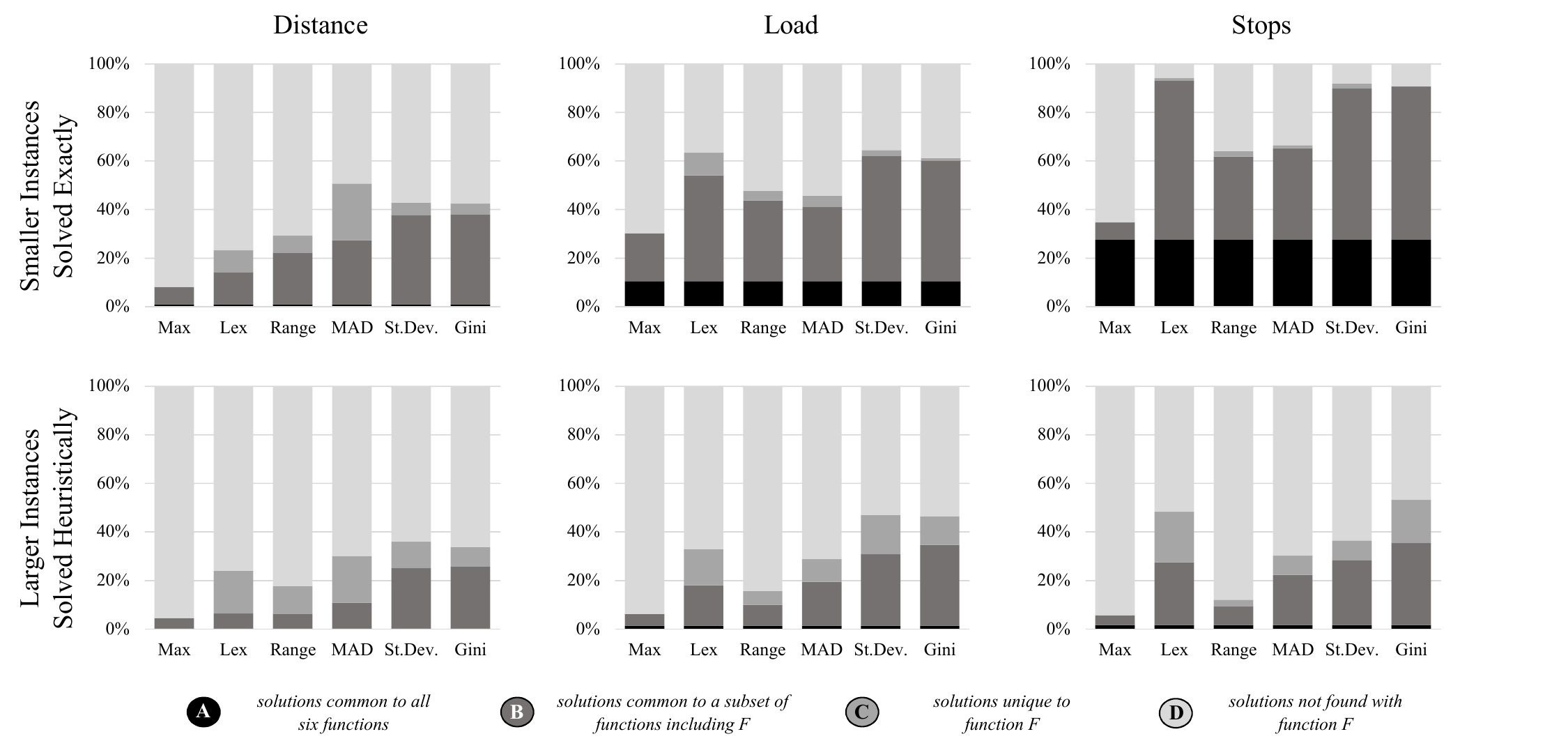}
	\caption{Overlap between non-dominated solution sets found with different equity functions balancing the same resource. The stacked bars indicate in relative terms the size of each category as a share of all the non-dominated solutions found for a given workload resource on all examined instances, i.e.~100\% corresponds to the union of all non-dominated solution sets identified with all equity functions for a given resource.}
	\label{overlap}
\end{figure}
\vfill
\end{landscape}

%\begin{landscape}
%\thispagestyle{empty}
%\begin{figure}
%	\centering
%		\includegraphics[width=1.5\textwidth]{figures/overlap.pdf}
%		\includegraphics[width=1.5\textwidth]{figures/overlap_color.pdf}
%	\caption{Overlap between non-dominated solution sets found with different equity functions balancing the same resource. The stacked bars indicate in relative terms the size of each category as a share of all the non-dominated solutions found for a given workload resource on all examined instances, i.e.~100\% corresponds to the union of all non-dominated solution sets identified with all equity functions for a given resource.}
%	\label{overlap}
%\end{figure}
%\end{landscape}

However, Figure \ref{overlap} also illustrates a certain paradox: if we consider the combined Pareto sets of each function in isolation (i.e.~categories A, B, and C only), then a significant majority of the identified solutions are shared with at least one other function. This would suggest that regardless of which equity function is chosen, it will find largely the same solutions as those found by the other functions. Yet looking at the complete set of all solutions for each model (i.e.~including category D), it is clear that the solutions \textit{not} found by our chosen function (category D) constitute the majority of all potentially ``well-balanced'' solutions, especially on the larger instances. This means that although no single function is particularly unique in terms of the solutions it finds, the fact that only \textit{one} equity function can be chosen means that any model will disregard a majority of solutions which would still be considered well-balanced otherwise.

%% file: tables/sample-sets.tex
\begin{table}[]
\centering
\scriptsize
\resizebox{\textwidth}{!}{%
\begin{tabular}{@{}ccccccccccccc@{}}
\toprule

\textbf{Sol.~Nr.} & \textbf{Cost} & \textbf{Max} & \textbf{Lex} & \textbf{Range} & \textbf{MAD} & \textbf{St.Dev.} & \textbf{Gini} & \textbf{Tour 1} & \textbf{Tour 2} & \textbf{Tour 3} & \textbf{Tour 4} & \textbf{Tour 5} \\ \midrule
1                 & 6654          & 279          & 1            & 186            & 52.08        & 67.89            & 0.1495        & 279             & 275             & 239             & 230             & 93              \\
2                 & 6715          & -            & 2            & -              & -            & 66.93            & 0.1430        & 279             & 257             & 257             & 230             & 93              \\
3                 & 6739          & 275          & -            & 182            & -            & -                & -             & 275             & 273             & 262             & 213             & 93              \\
4                 & 6739          & -            & 3            & -              & -            & 66.27            & 0.1387        & 275             & 262             & 247             & 239             & 93              \\
5                 & 6755          & -            & -            & 148            & 38.56        & 50.64            & 0.1197        & 279             & 257             & 230             & 219             & 131             \\
6                 & 6761          & 257          & 4            & 126            & 36.88        & 47.39            & 0.1007        & 257             & 255             & 247             & 226             & 131             \\
7                 & 6784          & -            & -            & 125            & -            & -                & -             & 279             & 275             & 239             & 169             & 154             \\
8                 & 6793          & -            & 5            & -              & -            & 46.97            & 0.0978        & 257             & 251             & 247             & 230             & 131             \\
9                 & 6839          & -            & -            & 110            & -            & 46.27            & -             & 279             & 275             & 213             & 180             & 169             \\
10                & 6862          & -            & -            & 106            & -            & -                & -             & 275             & 273             & 213             & 186             & 169             \\
11                & 6862          & -            & -            & -              & 36.56        & 39.55            & -             & 275             & 247             & 239             & 186             & 169             \\
12                & 6877          & -            & 6            & -              & -            & -                & 0.0957        & 257             & 248             & 247             & 233             & 131             \\
13                & 6880          & -            & -            & -              & 30.24        & 35.82            & 0.0885        & 275             & 247             & 213             & 212             & 169             \\
14                & 6900          & -            & 7            & -              & -            & -                & -             & 257             & 247             & 247             & 234             & 131             \\
15                & 6904          & -            & -            & 104            & -            & -                & -             & 273             & 257             & 231             & 186             & 169             \\
16                & 6922          & -            & 8            & 88             & 26.16        & 31.09            & 0.0756        & 257             & 247             & 231             & 212             & 169             \\
17                & 6941          & -            & 9            & -              & 25.76        & 30.97            & 0.0753        & 257             & 247             & 230             & 213             & 169             \\
18                & 6952          & 256          & 10           & -              & -            & -                & -             & 256             & 255             & 248             & 226             & 131             \\
19                & 6970          & 252          & 11           & 83             & -            & 30.41            & 0.0720        & 252             & 247             & 236             & 212             & 169             \\
20                & 6989          & 248          & 12           & 79             & -            & 30.09            & 0.0692        & 248             & 247             & 240             & 212             & 169             \\
21                & 6990          & -            & -            & 71             & -            & -                & -             & 257             & 247             & 239             & 187             & 186             \\
22                & 6991          & 247          & 13           & -              & -            & 29.99            & 0.0677        & 247             & 245             & 243             & 212             & 169             \\
23                & 7008          & -            & -            & 70             & -            & 25.49            & 0.0627        & 257             & 247             & 213             & 212             & 187             \\
24                & 7008          & -            & -            & -              & 22.16        & -                & -             & 248             & 247             & 230             & 222             & 169             \\
25                & 7038          & -            & 14           & -              & 21.68        & -                & 0.0581        & 247             & 236             & 234             & 230             & 169             \\
26                & 7057          & -            & -            & 62             & -            & -                & -             & 248             & 248             & 247             & 187             & 186             \\
27                & 7075          & -            & -            & 61             & 19.44        & 22.89            & 0.0563        & 248             & 247             & 222             & 212             & 187             \\
28                & 7105          & -            & 15           & 60             & 18.96        & 21.37            & 0.0516        & 247             & 236             & 234             & 212             & 187             \\
29                & 7151          & -            & -            & 58             & -            & -                & -             & 248             & 247             & 240             & 191             & 190             \\
30                & 7153          & -            & -            & 57             & -            & -                & -             & 247             & 245             & 243             & 191             & 190             \\
31                & 7166          & -            & 16           & -              & 16.56        & 20.33            & 0.0487        & 247             & 234             & 230             & 218             & 187             \\
32                & 7171          & -            & -            & 44             & -            & 20.03            & 0.0466        & 248             & 247             & 212             & 205             & 204             \\
33                & 7199          & -            & 17           & -              & 14.48        & 19.66            & 0.0441        & 247             & 229             & 227             & 226             & 187             \\
34                & 7201          & -            & -            & 42             & 13.84        & 15.28            & 0.0380        & 247             & 234             & 218             & 212             & 205             \\
35                & 7234          & -            & 18           & -              & 12.96        & -                & -             & 247             & 229             & 226             & 223             & 191             \\
36                & 7241          & -            & 19           & -              & -            & 15.17            & 0.0376        & 247             & 229             & 226             & 210             & 204             \\
37                & 7251          & 245          & 20           & -              & -            & -                & -             & 245             & 243             & 232             & 227             & 169             \\
38                & 7253          & -            & -            & 41             & -            & 15.07            & 0.0369        & 247             & 229             & 226             & 208             & 206             \\
39                & 7300          & -            & 21           & -              & -            & -                & -             & 245             & 236             & 234             & 232             & 169             \\
40                & 7301          & 239          & 22           & -              & -            & -                & -             & 239             & 238             & 236             & 234             & 169             \\
41                & 7351          & -            & -            & -              & 12.24        & 14.63            & 0.0366        & 247             & 230             & 222             & 212             & 205             \\
42                & 7365          & 236          & 23           & -              & -            & -                & -             & 236             & 234             & 232             & 227             & 187             \\
43                & 7365          & -            & -            & 39             & -            & -                & -             & 247             & 234             & 218             & 209             & 208             \\
44                & 7384          & -            & -            & -              & 12.16        & 14.55            & 0.0355        & 247             & 227             & 226             & 210             & 206             \\
45                & 7406          & -            & -            & -              & 10.64        & -                & -             & 247             & 226             & 223             & 220             & 200             \\
46                & 7424          & -            & -            & 37             & -            & 13.99            & -             & 247             & 230             & 219             & 210             & 210             \\
47                & 7424          & -            & -            & -              & -            & -                & 0.0330        & 257             & 217             & 216             & 213             & 213             \\
48                & 7431          & -            & -            & -              & -            & 13.59            & 0.0323        & 247             & 226             & 223             & 210             & 210             \\
49                & 7455          & -            & -            & -              & 9.52         & -                & 0.0319        & 247             & 223             & 223             & 218             & 205             \\
50                & 7460          & -            & 24           & -              & -            & -                & -             & 236             & 234             & 229             & 227             & 190             \\
51                & 7461          & 234          & 25           & 29             & 9.36         & 10.65            & 0.0258        & 234             & 232             & 227             & 218             & 205             \\
52                & 7497          & -            & -            & 28             & -            & -                & -             & 236             & 234             & 227             & 211             & 208             \\
53                & 7500          & -            & -            & -              & -            & 10.53            & -             & 236             & 234             & 221             & 217             & 208             \\
54                & 7521          & -            & -            & 23             & -            & 9.74             & 0.0229        & 236             & 234             & 217             & 216             & 213             \\
55                & 7550          & 230          & 26           & 18             & 6.16         & 6.73             & 0.0161        & 230             & 228             & 227             & 219             & 212             \\
56                & 7555          & 229          & 27           & 17             & 5.36         & 6.24             & 0.0147        & 229             & 228             & 226             & 221             & 212             \\
57                & 7557          & 228          & 28           & 16             & 4.56         & 5.84             & 0.0129        & 228             & 227             & 226             & 223             & 212             \\
58                & 7573          & -            & -            & 11             & 4.24         & 4.49             & 0.0108        & 230             & 227             & 221             & 219             & 219             \\
59                & 7578          & -            & -            & 10             & 3.44         & 3.71             & 0.0090        & 229             & 226             & 221             & 221             & 219             \\
60                & 7580          & 227          & 29           & 8              & 2.64         & 2.99             & 0.0075        & 227             & 226             & 223             & 221             & 219             \\
61                & 7689          & -            & 30           & 6              & 1.84         & 2.14             & 0.0050        & 227             & 224             & 222             & 222             & 221             \\
62                & 7884          & -            & 31           & 5              & 1.52         & 1.94             & 0.0039        & 227             & 223             & 222             & 222             & 222             \\
63                & 7909          & 226          & 32           & -              & -            & -                & -             & 226             & 224             & 224             & 221             & 221             \\
64                & 8062          & -            & 33           & -              & 1.44         & 1.72             & -             & 226             & 224             & 223             & 222             & 221             \\
65                & 8089          & 224          & 34           & 2              & 0.64         & 0.75             & 0.0018        & 224             & 224             & 223             & 223             & 222             \\
66                & 8274          & -            & 35           & 1              & 0.32         & 0.4              & 0.0007        & 224             & 223             & 223             & 223             & 223             \\ \bottomrule
\end{tabular}%
}
\caption{Pareto-optimal solution sets for a sample instance balancing tour loads. A dash indicates that the solution is not Pareto-optimal for the respective equity function.  The lexicographic objective is replaced with a ranking, as the workload vectors are given on the right.}
\label{sample-sets}
\end{table}

%% file: sections/4-5-agreement.tex
\subsection{Agreement Between Different Models}
\label{4-5}

Based on our analysis thus far, we have established that the workload resource has a notable impact on the trade-off structure between cost and balance, and that even for the same resource, no single equity function captures a majority of the potentially well-balanced allocations of that resource. However, all the presented resources and functions are intended to quantify minor variations of the same -- to an extent subjective -- concept. It therefore seems intuitive that even if the models do not produce \textit{exactly} the same solutions, there might be at least \textit{some} degree of agreement between them as to which VRP solutions are ``well-balanced'', and which are not. In this section we therefore examine to what extent the solutions optimizing one balance objective are still of high quality when considered for another objective. 

We quantify the degree of agreement as follows: For each instance, the Pareto-optimal or best known reference set is given for each of the 18 alternative balance objectives. 
For each objective, its Pareto set is re-evaluated according to each of the other 17 objectives, discarding any solutions which become dominated. We then compare these 17 re-evaluated Pareto sets to the corresponding optimal Pareto set of each objective, using the well-known hypervolume indicator as a measure of overall solution quality (we refer the reader to \cite{zitzler2008} for a detailed discussion of quality indicators for multi-objective optimization). The hypervolume represents the volume of the objective space dominated by a solution set, so the optimal Pareto set for each objective corresponds to 100\% of the attainable hypervolume. We note that this indicator is sensitive to the choice of the reference point delimiting the objective space -- we therefore use the nadir point of each objective's optimal Pareto set since this is the strictest valid reference point. 

Tables \ref{agreement1} and \ref{agreement2} report the results of the above analysis on the set of smaller instances solved to optimality, and the larger ones solved heuristically. The values in the tables represent the attained hypervolume as a percentage of the optimal or best known hypervolume, averaged over all the instances in the respective set. Values of 90\% or more are shaded and printed in bold.

\input{tables/agreement1}
\input{tables/agreement2}
%\onehalfspacing

The most immediate observation we can make is that solutions balancing one resource are generally not of good quality in terms of balancing the other potential resources. Nearly all the inter-resource comparisons result in hypervolume values well below 90\%, and the values worsen further on the larger instances. In Section \ref{4-3} we identified noticeable differences in the trade-off structure when the resource is changed, but those differences were considerably smaller than those observed here. This suggests that although well-balanced solutions can be found at low cost for all the examined resources, those solutions tend to differ quite dramatically in structure despite sharing very similar cost totals.

Compared to the effect of changing the workload resource, using a different equity function has only limited impact. We can see that most of the intra-resource comparisons in Tables \ref{agreement1} and \ref{agreement2} result in values well above 90\%, in stark contrast to the inter-resource comparisons. Although this reflects the trends observed in Section \ref{4-3} with regard to similar trade-off structures, it is not a trivial outcome given that in Section \ref{4-4} the degree of strict overlap between any two functions was found to be much lower. This implies that for a given resource, the Pareto-optimal solutions generated by different equity functions are seldom identical, but their quality is often similar. In other words, high-quality solutions balancing the same resource tend to be of similar quality regardless of the equity function, but not vice versa. 

%The implication is that although different functions generally do not identify exactly the same solutions as Pareto-optimal given the same resource, the solutions they do identify tend to be much more similar than when the resource is changed and not the function. In other words, solutions balancing the same resource tend to be similar regardless of the function, but not vice versa.

Tables \ref{agreement1} and \ref{agreement2} also illustrate two relevant exceptions to the above trends. Looking at the comparisons within the distance resource, we can see that the two monotonic functions (max, lexicographic) are not particularly good at approximating the types of solutions found with non-monotonic functions, and vice versa. Although we accounted for the TSP-optimality and workload consistency issues that arise when balancing a variable-sum resource with non-monotonic equity functions, it is clear that significant differences remain. The other important observation is that simpler functions like max and range are less effective at providing good approximations of the other Pareto sets, even for the same resource. This is most likely a direct consequence of the lower cardinality of the Pareto sets identified with these functions, which reduces the probability of finding solutions which are also of high quality for other functions. This also makes it easier for the more sophisticated functions to better approximate the small Pareto sets found with the max and range functions -- this is particularly visible for the stops resource on the smaller instances.

%% file: tables/agreement1.tex
%CCCCCCCCCCCCCCCCCC
%cccccccccccccccccc

%\newgeometry{margin=2cm}
\begin{landscape}
\vspace*{\fill}
\thispagestyle{empty}

\begin{table}[htb]
\centering
\resizebox{\linewidth}{!}{%
\begin{tabular}{cclCCCCCCCCCCCCCCCCCC}
                                               &                            &         & \multicolumn{18}{c}{\textbf{Re-evaluation Objective $B$}}                                                                                                                                                                                                                                                                                                                                                                                                                                                                                                                                                                                                                                                                                                                                 \\
                                               &                            &         & \multicolumn{6}{c|}{Distance}                                                                                                                                                                                                                                    & \multicolumn{6}{c|}{Load}                                                                                                                                                                                                                                          & \multicolumn{6}{c}{Stops}                                                                                                                                                                                                                     \\
                                               &                            &         & Max                                  & Lex                                  & Range                                 & MAD                                   & St.Dev.                               & \multicolumn{1}{c|}{Gini}                                  & Max                                   & Lex                                   & Range                                 & MAD                                   & St.Dev.                               & \multicolumn{1}{c|}{Gini}                                  & Max                                   & Lex                                   & Range                                 & MAD                                   & St.Dev.                               & Gini                                  \\
                                               &                            & Max     & \cellcolor[HTML]{EFEFEF}\textbf{100} & \cellcolor[HTML]{EFEFEF}\textbf{99.0} & 67.4                                  & 65.1                                  & 65.0                                  & \multicolumn{1}{c|}{73.6}                                  & 26.9                                  & 25.1                                  & 44.7                                  & 44.3                                  & 45.9                                  & \multicolumn{1}{c|}{40.6}                                  & 89.4                                  & 18.6                                  & 67.5                                  & 42.0                                  & 41.7                                  & 36.2                                  \\
                                               &                            & Lex     & \cellcolor[HTML]{EFEFEF}\textbf{100} & \cellcolor[HTML]{EFEFEF}\textbf{100} & 82.5                                  & 81.5                                  & 81.7                                  & \multicolumn{1}{c|}{86.7}                                  & 41.0                                  & 39.5                                  & 61.1                                  & 59.4                                  & 61.4                                  & \multicolumn{1}{c|}{56.5}                                  & \cellcolor[HTML]{EFEFEF}\textbf{93.3} & 29.4                                  & 84.8                                  & 69.1                                  & 69.8                                  & 64.2                                  \\
                                               &                            & Range   & 81.0                                 & 84.9                                 & \cellcolor[HTML]{EFEFEF}\textbf{100}  & \cellcolor[HTML]{EFEFEF}\textbf{93.5} & \cellcolor[HTML]{EFEFEF}\textbf{98.0} & \multicolumn{1}{c|}{\cellcolor[HTML]{EFEFEF}\textbf{98.0}} & 48.5                                  & 47.3                                  & 70.2                                  & 69.2                                  & 70.9                                  & \multicolumn{1}{c|}{66.4}                                  & \cellcolor[HTML]{EFEFEF}\textbf{93.5} & 37.2                                  & \cellcolor[HTML]{EFEFEF}\textbf{91.4} & 69.9                                  & 76.5                                  & 69.8                                  \\
                                               &                            & MAD     & 66.6                                 & 71.7                                 & \cellcolor[HTML]{EFEFEF}\textbf{95.1} & \cellcolor[HTML]{EFEFEF}\textbf{100}  & \cellcolor[HTML]{EFEFEF}\textbf{98.7} & \multicolumn{1}{c|}{\cellcolor[HTML]{EFEFEF}\textbf{98.6}} & 51.3                                  & 50.0                                  & 71.4                                  & 72.4                                  & 72.9                                  & \multicolumn{1}{c|}{68.8}                                  & \cellcolor[HTML]{EFEFEF}\textbf{95.3} & 45.5                                  & \cellcolor[HTML]{EFEFEF}\textbf{92.4} & 80.5                                  & 82.8                                  & 79.5                                  \\
                                               &                            & St.Dev. & 73.7                                 & 78.0                                 & \cellcolor[HTML]{EFEFEF}\textbf{98.8} & \cellcolor[HTML]{EFEFEF}\textbf{98.4} & \cellcolor[HTML]{EFEFEF}\textbf{100}  & \multicolumn{1}{c|}{\cellcolor[HTML]{EFEFEF}\textbf{99.6}} & 53.6                                  & 52.5                                  & 73.1                                  & 73.0                                  & 74.2                                  & \multicolumn{1}{c|}{70.2}                                  & \cellcolor[HTML]{EFEFEF}\textbf{95.5} & 53.2                                  & \cellcolor[HTML]{EFEFEF}\textbf{94.5} & 80.0                                  & 84.0                                  & 80.1                                  \\
                                               & \multirow{-6}{*}{Distance} & Gini    & 77.5                                 & 80.3                                 & \cellcolor[HTML]{EFEFEF}\textbf{98.8} & \cellcolor[HTML]{EFEFEF}\textbf{98.4} & \cellcolor[HTML]{EFEFEF}\textbf{99.6} & \multicolumn{1}{c|}{\cellcolor[HTML]{EFEFEF}\textbf{100}}  & 53.9                                  & 52.9                                  & 72.8                                  & 72.7                                  & 73.8                                  & \multicolumn{1}{c|}{69.9}                                  & \cellcolor[HTML]{EFEFEF}\textbf{95.4} & 52.8                                  & \cellcolor[HTML]{EFEFEF}\textbf{94.3} & 80.1                                  & 84.3                                  & 80.7                                  \\ \hhline{~--|------|------|------}
                                               &                            & Max     & 40.0                                 & 37.7                                 & 55.5                                  & 56.4                                  & 56.1                                  & \multicolumn{1}{c|}{61.4}                                  & \cellcolor[HTML]{EFEFEF}\textbf{100}  & \cellcolor[HTML]{EFEFEF}\textbf{100}  & \cellcolor[HTML]{EFEFEF}\textbf{95.3} & \cellcolor[HTML]{EFEFEF}\textbf{95.0} & \cellcolor[HTML]{EFEFEF}\textbf{94.9} & \multicolumn{1}{c|}{\cellcolor[HTML]{EFEFEF}\textbf{96.0}} & \cellcolor[HTML]{EFEFEF}\textbf{91.7} & 45.8                                  & 81.7                                  & 68.4                                  & 66.9                                  & 62.2                                  \\
                                               &                            & Lex     & 45.4                                 & 44.4                                 & 60.9                                  & 61.9                                  & 61.3                                  & \multicolumn{1}{c|}{66.7}                                  & \cellcolor[HTML]{EFEFEF}\textbf{100}  & \cellcolor[HTML]{EFEFEF}\textbf{100}  & \cellcolor[HTML]{EFEFEF}\textbf{97.8} & \cellcolor[HTML]{EFEFEF}\textbf{97.6} & \cellcolor[HTML]{EFEFEF}\textbf{97.6} & \multicolumn{1}{c|}{\cellcolor[HTML]{EFEFEF}\textbf{98.4}} & \cellcolor[HTML]{EFEFEF}\textbf{93.2} & 51.7                                  & 88.8                                  & 79.0                                  & 78.2                                  & 74.1                                  \\
                                               &                            & Range   & 43.7                                 & 41.5                                 & 64.7                                  & 65.3                                  & 65.3                                  & \multicolumn{1}{c|}{69.5}                                  & \cellcolor[HTML]{EFEFEF}\textbf{98.1} & \cellcolor[HTML]{EFEFEF}\textbf{98.2} & \cellcolor[HTML]{EFEFEF}\textbf{100}  & \cellcolor[HTML]{EFEFEF}\textbf{99.0} & \cellcolor[HTML]{EFEFEF}\textbf{99.6} & \multicolumn{1}{c|}{\cellcolor[HTML]{EFEFEF}\textbf{99.5}} & \cellcolor[HTML]{EFEFEF}\textbf{93.0} & 63.4                                  & \cellcolor[HTML]{EFEFEF}\textbf{94.1} & 83.8                                  & 86.6                                  & 84.3                                  \\
                                               &                            & MAD     & 40.0                                 & 39.9                                 & 63.9                                  & 64.0                                  & 64.5                                  & \multicolumn{1}{c|}{68.9}                                  & \cellcolor[HTML]{EFEFEF}\textbf{95.4} & \cellcolor[HTML]{EFEFEF}\textbf{95.6} & \cellcolor[HTML]{EFEFEF}\textbf{98.7} & \cellcolor[HTML]{EFEFEF}\textbf{100}  & \cellcolor[HTML]{EFEFEF}\textbf{99.6} & \multicolumn{1}{c|}{\cellcolor[HTML]{EFEFEF}\textbf{99.2}} & \cellcolor[HTML]{EFEFEF}\textbf{94.2} & 61.4                                  & \cellcolor[HTML]{EFEFEF}\textbf{93.5} & 83.9                                  & 86.1                                  & 83.4                                  \\
                                               &                            & St.Dev. & 45.2                                 & 43.5                                 & 66.0                                  & 66.3                                  & 66.2                                  & \multicolumn{1}{c|}{70.6}                                  & \cellcolor[HTML]{EFEFEF}\textbf{98.2} & \cellcolor[HTML]{EFEFEF}\textbf{98.3} & \cellcolor[HTML]{EFEFEF}\textbf{99.9} & \cellcolor[HTML]{EFEFEF}\textbf{99.9} & \cellcolor[HTML]{EFEFEF}\textbf{100}  & \multicolumn{1}{c|}{\cellcolor[HTML]{EFEFEF}\textbf{100}}  & \cellcolor[HTML]{EFEFEF}\textbf{95.2} & 65.7                                  & \cellcolor[HTML]{EFEFEF}\textbf{95.4} & 86.9                                  & 89.5                                  & 87.1                                  \\
                                               & \multirow{-6}{*}{Load}     & Gini    & 44.7                                 & 42.8                                 & 65.5                                  & 66.3                                  & 66.1                                  & \multicolumn{1}{c|}{70.4}                                  & \cellcolor[HTML]{EFEFEF}\textbf{98.4} & \cellcolor[HTML]{EFEFEF}\textbf{98.5} & \cellcolor[HTML]{EFEFEF}\textbf{99.8} & \cellcolor[HTML]{EFEFEF}\textbf{99.9} & \cellcolor[HTML]{EFEFEF}\textbf{99.9} & \multicolumn{1}{c|}{\cellcolor[HTML]{EFEFEF}\textbf{100}}  & \cellcolor[HTML]{EFEFEF}\textbf{95.0} & 63.0                                  & \cellcolor[HTML]{EFEFEF}\textbf{94.7} & 86.0                                  & 88.3                                  & 85.9                                  \\ \hhline{~--|------|------|------}
                                               &                            & Max     & 24.1                                 & 23.7                                 & 42.4                                  & 43.4                                  & 44.6                                  & \multicolumn{1}{c|}{46.5}                                  & 39.4                                  & 37.8                                  & 58.4                                  & 56.0                                  & 58.4                                  & \multicolumn{1}{c|}{54.0}                                  & \cellcolor[HTML]{EFEFEF}\textbf{100}  & \cellcolor[HTML]{EFEFEF}\textbf{98.3} & 66.7                                  & 24.8                                  & 26.7                                  & 28.6                                  \\
                                               &                            & Lex     & 34.2                                 & 32.4                                 & 50.6                                  & 50.1                                  & 52.0                                  & \multicolumn{1}{c|}{53.9}                                  & 45.5                                  & 43.9                                  & 67.3                                  & 64.9                                  & 67.5                                  & \multicolumn{1}{c|}{62.4}                                  & \cellcolor[HTML]{EFEFEF}\textbf{100}  & \cellcolor[HTML]{EFEFEF}\textbf{100}  & \cellcolor[HTML]{EFEFEF}\textbf{99.9} & \cellcolor[HTML]{EFEFEF}\textbf{99.9} & \cellcolor[HTML]{EFEFEF}\textbf{99.8} & \cellcolor[HTML]{EFEFEF}\textbf{100}  \\
                                               &                            & Range   & 30.5                                 & 29.1                                 & 49.4                                  & 47.5                                  & 50.0                                  & \multicolumn{1}{c|}{51.8}                                  & 42.2                                  & 40.6                                  & 65.6                                  & 62.2                                  & 65.3                                  & \multicolumn{1}{c|}{59.8}                                  & \cellcolor[HTML]{EFEFEF}\textbf{99.6} & \cellcolor[HTML]{EFEFEF}\textbf{97.9} & \cellcolor[HTML]{EFEFEF}\textbf{100}  & 63.9                                  & 77.6                                  & 73.7                                  \\
                                               &                            & MAD     & 32.4                                 & 30.8                                 & 49.8                                  & 49.7                                  & 51.4                                  & \multicolumn{1}{c|}{53.3}                                  & 43.7                                  & 42.1                                  & 65.5                                  & 63.8                                  & 66.2                                  & \multicolumn{1}{c|}{61.1}                                  & \cellcolor[HTML]{EFEFEF}\textbf{96.8} & \cellcolor[HTML]{EFEFEF}\textbf{93.1} & \cellcolor[HTML]{EFEFEF}\textbf{93.4} & \cellcolor[HTML]{EFEFEF}\textbf{100}  & \cellcolor[HTML]{EFEFEF}\textbf{94.0} & \cellcolor[HTML]{EFEFEF}\textbf{94.9} \\
                                               &                            & St.Dev. & 34.0                                 & 32.3                                 & 50.7                                  & 50.2                                  & 52.0                                  & \multicolumn{1}{c|}{53.9}                                  & 45.5                                  & 43.8                                  & 67.3                                  & 65.0                                  & 67.5                                  & \multicolumn{1}{c|}{62.4}                                  & \cellcolor[HTML]{EFEFEF}\textbf{99.7} & \cellcolor[HTML]{EFEFEF}\textbf{99.3} & \cellcolor[HTML]{EFEFEF}\textbf{99.9} & \cellcolor[HTML]{EFEFEF}\textbf{100}  & \cellcolor[HTML]{EFEFEF}\textbf{100}  & \cellcolor[HTML]{EFEFEF}\textbf{100}  \\
\multirow{-18}{*}{\rotatebox[origin=c]{90}{\textbf{Original Objective $A$}}} & \multirow{-6}{*}{Stops}    & Gini    & 34.0                                 & 32.3                                 & 50.6                                  & 50.1                                  & 52.0                                  & \multicolumn{1}{c|}{53.9}                                  & 45.5                                  & 43.8                                  & 67.3                                  & 64.9                                  & 67.5                                  & \multicolumn{1}{c|}{62.4}                                  & \cellcolor[HTML]{EFEFEF}\textbf{99.7} & \cellcolor[HTML]{EFEFEF}\textbf{99.4} & \cellcolor[HTML]{EFEFEF}\textbf{99.9} & \cellcolor[HTML]{EFEFEF}\textbf{100}  & \cellcolor[HTML]{EFEFEF}\textbf{99.9} & \cellcolor[HTML]{EFEFEF}\textbf{100}
\end{tabular}%
}
\caption{Hypervolume values attained by Pareto-optimal solution sets identified by an objective $A$ and re-evaluated according to an alternative objective~$B$. Values represent the attained hypervolume as a percentage of the optimal hypervolume, averaged over the set of smaller instances solved to optimality. Values of 90\% or more are shaded and typeset in bold.}
\label{agreement1}
\end{table}

\vfill
\end{landscape}

%\restoregeometry

%% file: tables/agreement2.tex
%\newgeometry{margin=2cm}

\begin{landscape}
\vspace*{\fill}
\thispagestyle{empty}
\begin{table}[htb]
\centering
\resizebox{\linewidth}{!}{%
\begin{tabular}{cclCCCCCCCCCCCCCCCCCC}
                                                     &                            &         & \multicolumn{18}{c}{\textbf{Re-evaluation Objective $B$}}                                                                                                                                                                                                                                                                                                                                                                                                                                                                                                                                                                                                                                                                                                                                               \\
                                                     &                            &         & \multicolumn{6}{c|}{Distance}                                                                                                                                                                                                                                            & \multicolumn{6}{c|}{Load}                                                                                                                                                                                                                                                & \multicolumn{6}{c}{Stops}                                                                                                                                                                                                                           \\
                                                     &                            &         & Max                                    & Lex                                    & Range                                  & MAD                                    & St.Dev.                                & \multicolumn{1}{c|}{Gini}                                   & Max                                    & Lex                                    & Range                                  & MAD                                    & St.Dev.                                & \multicolumn{1}{c|}{Gini}                                   & Max                                    & Lex                                    & Range                                  & MAD                                    & St.Dev.                                & Gini                                   \\
                                                     &                            & Max     & \cellcolor[HTML]{EFEFEF}\textbf{100.0} & \cellcolor[HTML]{EFEFEF}\textbf{98.6}  & 51.0                                   & 41.7                                   & 43.2                                   & \multicolumn{1}{c|}{49.8}                                   & 35.6                                   & 6.9                                    & 36.6                                   & 25.9                                   & 29.7                                   & \multicolumn{1}{c|}{23.1}                                   & 72.4                                   & 27.8                                   & 51.1                                   & 23.8                                   & 24.9                                   & 22.8                                   \\
                                                     &                            & Lex     & \cellcolor[HTML]{EFEFEF}\textbf{100.0} & \cellcolor[HTML]{EFEFEF}\textbf{100.0} & 68.9                                   & 67.9                                   & 67.6                                   & \multicolumn{1}{c|}{72.7}                                   & 39.5                                   & 11.5                                   & 46.7                                   & 37.0                                   & 41.1                                   & \multicolumn{1}{c|}{33.0}                                   & 79.6                                   & 40.1                                   & 61.4                                   & 36.8                                   & 38.0                                   & 35.5                                   \\
                                                     &                            & Range   & 81.8                                   & 82.2                                   & \cellcolor[HTML]{EFEFEF}\textbf{100.0} & 89.0                                   & \cellcolor[HTML]{EFEFEF}\textbf{94.9}  & \multicolumn{1}{c|}{\cellcolor[HTML]{EFEFEF}\textbf{95.9}}  & 39.4                                   & 11.2                                   & 55.2                                   & 42.4                                   & 49.7                                   & \multicolumn{1}{c|}{38.2}                                   & 75.0                                   & 32.4                                   & 66.6                                   & 40.7                                   & 43.9                                   & 38.8                                   \\
                                                     &                            & MAD     & 60.9                                   & 64.3                                   & \cellcolor[HTML]{EFEFEF}\textbf{90.6}  & \cellcolor[HTML]{EFEFEF}\textbf{100.0} & \cellcolor[HTML]{EFEFEF}\textbf{97.3}  & \multicolumn{1}{c|}{\cellcolor[HTML]{EFEFEF}\textbf{98.2}}  & 41.8                                   & 13.9                                   & 56.9                                   & 47.1                                   & 52.4                                   & \multicolumn{1}{c|}{42.7}                                   & 76.0                                   & 33.1                                   & 64.4                                   & 43.3                                   & 43.5                                   & 41.3                                   \\
                                                     &                            & St.Dev. & 69.2                                   & 71.6                                   & \cellcolor[HTML]{EFEFEF}\textbf{97.2}  & \cellcolor[HTML]{EFEFEF}\textbf{98.0}  & \cellcolor[HTML]{EFEFEF}\textbf{100.0} & \multicolumn{1}{c|}{\cellcolor[HTML]{EFEFEF}\textbf{99.5}}  & 43.8                                   & 16.4                                   & 60.7                                   & 48.4                                   & 55.4                                   & \multicolumn{1}{c|}{44.1}                                   & 77.3                                   & 35.2                                   & 70.0                                   & 47.1                                   & 49.3                                   & 45.4                                   \\
                                                     & \multirow{-6}{*}{Distance} & Gini    & 69.9                                   & 71.3                                   & \cellcolor[HTML]{EFEFEF}\textbf{97.0}  & \cellcolor[HTML]{EFEFEF}\textbf{98.3}  & \cellcolor[HTML]{EFEFEF}\textbf{99.3}  & \multicolumn{1}{c|}{\cellcolor[HTML]{EFEFEF}\textbf{100.0}} & 42.3                                   & 14.3                                   & 60.0                                   & 48.3                                   & 54.8                                   & \multicolumn{1}{c|}{43.9}                                   & 78.5                                   & 38.5                                   & 69.0                                   & 47.0                                   & 48.0                                   & 45.0                                   \\ \hhline{~--|------|------|------}
                                                     &                            & Max     & 32.2                                   & 26.0                                   & 25.9                                   & 25.0                                   & 24.9                                   & \multicolumn{1}{c|}{27.3}                                   & \cellcolor[HTML]{EFEFEF}\textbf{100.0} & \cellcolor[HTML]{EFEFEF}\textbf{97.1}  & 77.7                                   & 75.3                                   & 73.9                                   & \multicolumn{1}{c|}{77.1}                                   & 69.2                                   & 21.9                                   & 54.3                                   & 24.9                                   & 27.2                                   & 24.6                                   \\
                                                     &                            & Lex     & 40.7                                   & 33.8                                   & 32.8                                   & 32.0                                   & 31.9                                   & \multicolumn{1}{c|}{34.2}                                   & \cellcolor[HTML]{EFEFEF}\textbf{100.0} & \cellcolor[HTML]{EFEFEF}\textbf{100.0} & \cellcolor[HTML]{EFEFEF}\textbf{94.0}  & \cellcolor[HTML]{EFEFEF}\textbf{96.2}  & \cellcolor[HTML]{EFEFEF}\textbf{95.6}  & \multicolumn{1}{c|}{\cellcolor[HTML]{EFEFEF}\textbf{97.3}}  & 76.7                                   & 34.9                                   & 65.6                                   & 43.4                                   & 44.7                                   & 41.9                                   \\
                                                     &                            & Range   & 38.8                                   & 31.9                                   & 33.8                                   & 32.0                                   & 32.6                                   & \multicolumn{1}{c|}{34.0}                                   & \cellcolor[HTML]{EFEFEF}\textbf{91.0}  & 88.2                                   & \cellcolor[HTML]{EFEFEF}\textbf{100.0} & \cellcolor[HTML]{EFEFEF}\textbf{92.1}  & \cellcolor[HTML]{EFEFEF}\textbf{96.0}  & \multicolumn{1}{c|}{\cellcolor[HTML]{EFEFEF}\textbf{93.8}}  & 71.8                                   & 28.9                                   & 67.0                                   & 40.6                                   & 43.8                                   & 39.9                                   \\
                                                     &                            & MAD     & 39.7                                   & 32.9                                   & 33.7                                   & 32.3                                   & 32.7                                   & \multicolumn{1}{c|}{34.6}                                   & 87.7                                   & 83.9                                   & \cellcolor[HTML]{EFEFEF}\textbf{95.2}  & \cellcolor[HTML]{EFEFEF}\textbf{100.0} & \cellcolor[HTML]{EFEFEF}\textbf{98.2}  & \multicolumn{1}{c|}{\cellcolor[HTML]{EFEFEF}\textbf{98.4}}  & 75.6                                   & 33.9                                   & 67.5                                   & 46.1                                   & 47.7                                   & 44.8                                   \\
                                                     &                            & St.Dev. & 41.1                                   & 33.9                                   & 34.7                                   & 33.5                                   & 33.7                                   & \multicolumn{1}{c|}{35.4}                                   & \cellcolor[HTML]{EFEFEF}\textbf{91.2}  & 88.9                                   & \cellcolor[HTML]{EFEFEF}\textbf{99.4}  & \cellcolor[HTML]{EFEFEF}\textbf{99.2}  & \cellcolor[HTML]{EFEFEF}\textbf{100.0} & \multicolumn{1}{c|}{\cellcolor[HTML]{EFEFEF}\textbf{99.5}}  & 77.0                                   & 36.1                                   & 70.4                                   & 47.7                                   & 49.0                                   & 47.0                                   \\
                                                     & \multirow{-6}{*}{Load}     & Gini    & 39.4                                   & 32.6                                   & 34.4                                   & 33.1                                   & 33.3                                   & \multicolumn{1}{c|}{35.3}                                   & \cellcolor[HTML]{EFEFEF}\textbf{92.9}  & \cellcolor[HTML]{EFEFEF}\textbf{90.2}  & \cellcolor[HTML]{EFEFEF}\textbf{98.5}  & \cellcolor[HTML]{EFEFEF}\textbf{99.6}  & \cellcolor[HTML]{EFEFEF}\textbf{99.5}  & \multicolumn{1}{c|}{\cellcolor[HTML]{EFEFEF}\textbf{100.0}} & 77.6                                   & 36.9                                   & 69.5                                   & 47.1                                   & 49.3                                   & 46.5                                   \\ \hhline{~--|------|------|------}
                                                     &                            & Max     & 27.3                                   & 19.3                                   & 20.9                                   & 17.4                                   & 18.0                                   & \multicolumn{1}{c|}{18.9}                                   & 34.2                                   & 6.0                                    & 34.5                                   & 23.5                                   & 28.0                                   & \multicolumn{1}{c|}{21.0}                                   & \cellcolor[HTML]{EFEFEF}\textbf{100.0} & \cellcolor[HTML]{EFEFEF}\textbf{95.7}  & 74.0                                   & 47.0                                   & 50.0                                   & 53.3                                   \\
                                                     &                            & Lex     & 36.3                                   & 29.3                                   & 32.2                                   & 27.8                                   & 29.2                                   & \multicolumn{1}{c|}{29.7}                                   & 40.9                                   & 12.7                                   & 54.1                                   & 43.1                                   & 49.1                                   & \multicolumn{1}{c|}{39.0}                                   & \cellcolor[HTML]{EFEFEF}\textbf{100.0} & \cellcolor[HTML]{EFEFEF}\textbf{100.0} & \cellcolor[HTML]{EFEFEF}\textbf{97.8}  & \cellcolor[HTML]{EFEFEF}\textbf{97.3}  & \cellcolor[HTML]{EFEFEF}\textbf{98.2}  & \cellcolor[HTML]{EFEFEF}\textbf{98.2}  \\
                                                     &                            & Range   & 30.5                                   & 22.8                                   & 29.5                                   & 25.6                                   & 26.9                                   & \multicolumn{1}{c|}{27.8}                                   & 37.7                                   & 9.2                                    & 52.7                                   & 37.7                                   & 45.5                                   & \multicolumn{1}{c|}{34.3}                                   & \cellcolor[HTML]{EFEFEF}\textbf{96.8}  & \cellcolor[HTML]{EFEFEF}\textbf{90.2}  & \cellcolor[HTML]{EFEFEF}\textbf{100.0} & 70.9                                   & 81.5                                   & 77.7                                   \\
                                                     &                            & MAD     & 34.3                                   & 28.3                                   & 32.4                                   & 27.8                                   & 29.4                                   & \multicolumn{1}{c|}{30.0}                                   & 39.9                                   & 11.4                                   & 53.7                                   & 41.9                                   & 48.2                                   & \multicolumn{1}{c|}{38.3}                                   & \cellcolor[HTML]{EFEFEF}\textbf{91.7}  & 88.4                                   & \cellcolor[HTML]{EFEFEF}\textbf{92.4}  & \cellcolor[HTML]{EFEFEF}\textbf{100.0} & \cellcolor[HTML]{EFEFEF}\textbf{96.3}  & \cellcolor[HTML]{EFEFEF}\textbf{98.1}  \\
                                                     &                            & St.Dev. & 35.0                                   & 27.9                                   & 32.5                                   & 27.7                                   & 29.4                                   & \multicolumn{1}{c|}{29.7}                                   & 40.2                                   & 11.7                                   & 55.3                                   & 43.3                                   & 50.2                                   & \multicolumn{1}{c|}{39.6}                                   & \cellcolor[HTML]{EFEFEF}\textbf{96.9}  & \cellcolor[HTML]{EFEFEF}\textbf{92.2}  & \cellcolor[HTML]{EFEFEF}\textbf{98.6}  & \cellcolor[HTML]{EFEFEF}\textbf{97.9}  & \cellcolor[HTML]{EFEFEF}\textbf{100.0} & \cellcolor[HTML]{EFEFEF}\textbf{98.3}  \\
\multirow{-18}{*}{\rotatebox[origin=c]{90}{\textbf{Original Objective $A$}}} & \multirow{-6}{*}{Stops}    & Gini    & 35.7                                   & 28.8                                   & 32.8                                   & 28.4                                   & 29.7                                   & \multicolumn{1}{c|}{30.2}                                   & 40.6                                   & 12.3                                   & 55.4                                   & 44.3                                   & 50.6                                   & \multicolumn{1}{c|}{40.4}                                   & \cellcolor[HTML]{EFEFEF}\textbf{96.3}  & \cellcolor[HTML]{EFEFEF}\textbf{94.7}  & \cellcolor[HTML]{EFEFEF}\textbf{97.7}  & \cellcolor[HTML]{EFEFEF}\textbf{99.4}  & \cellcolor[HTML]{EFEFEF}\textbf{99.1}  & \cellcolor[HTML]{EFEFEF}\textbf{100.0}
\end{tabular}%
}
\caption{Hypervolume values attained by Pareto-optimal solution sets identified by an objective $A$ and re-evaluated according to an alternative objective~$B$. Values represent the attained hypervolume as a percentage of the optimal hypervolume, averaged over the set of larger instances solved heuristically. Values of 90\% or more are shaded and typeset in bold.}
\label{agreement2}
\end{table}
\vfill
\end{landscape}

%\restoregeometry

%% file: sections/4-6-similarity.tex
\subsection{Solution Similarity in the Decision Space}
\label{4-6}

We close our numerical study by considering the degree to which different Pareto-optimal VRP solutions for the same model share a common solution structure. Most multi-objective methods -- and particularly genetic algorithms -- implicitly assume that Pareto-optimal solutions share common characteristics, especially if those solutions have similar objective function values. This intuition has been a subject of research in the broader literature on fitness landscape analysis \citep{pitzer2012}, also for single-objective problems. In this last section we examine -- in our specific context of balanced VRPs -- whether there is a connection between solution similarity in the objective space and similarity in the decision space.

In the following, we quantify solution similarity as the percentage of edges common to a pair of solutions from the same Pareto set. For the instances used in our study, all vehicles are required and so all solutions of a given instance have the same number of edges, namely $n+K$. With $E(x_i)$ being the set of edges used in solution $x_i$, we define the similarity of a pair of solutions $i,j$ as:

\[ s(x_i, x_j) = \frac{|E(x_i) \cap E(x_j)|}{n + K} \]

\begin{figure}
	\centering
		\includegraphics[width=1.00\textwidth]{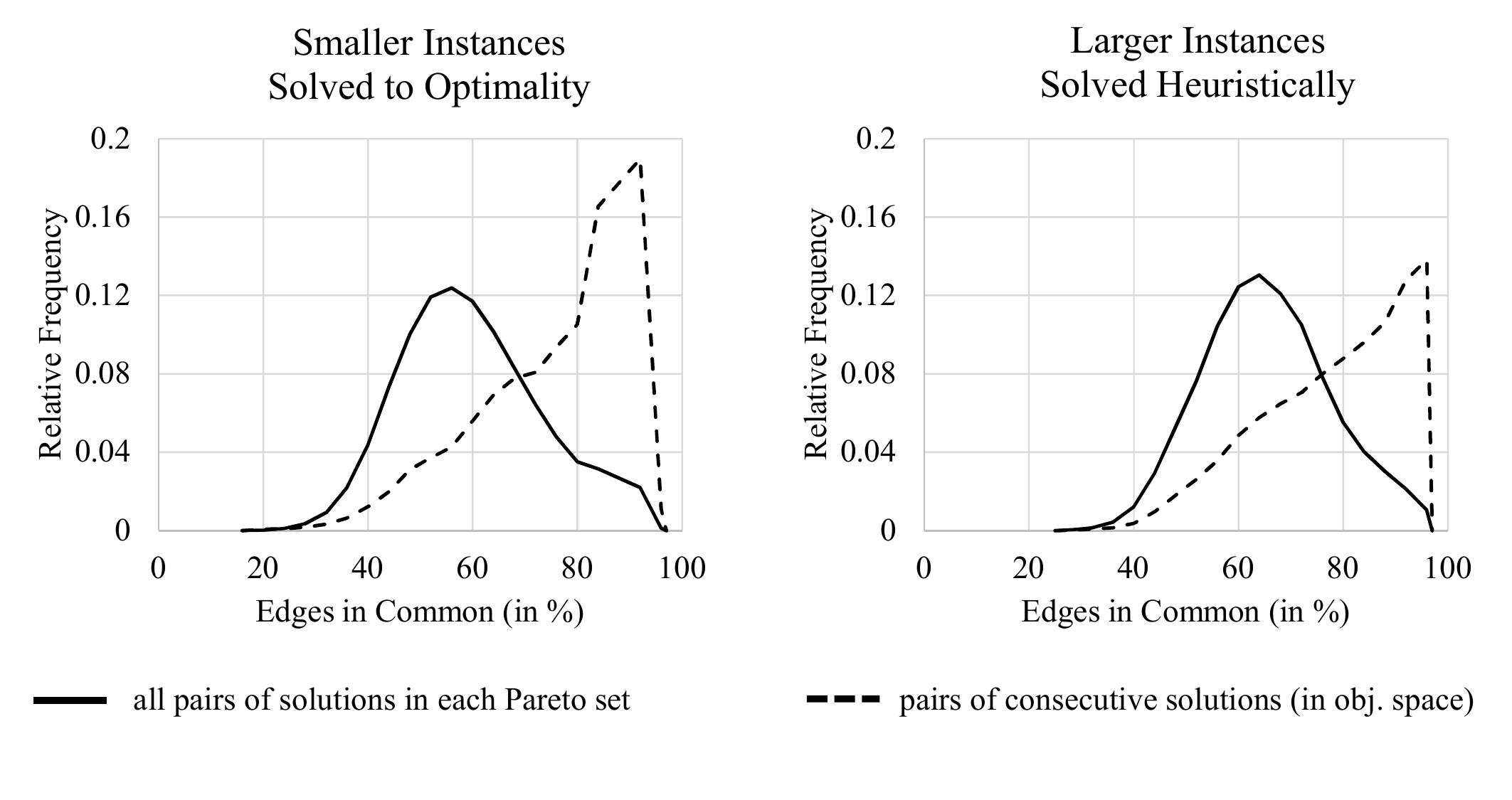}
	\caption{Frequency distributions of solution similarity between different pairs of balanced VRP solutions. Solution similarity is measured as the percentage of edges shared by a pair of solutions.}
	\label{similarity}
\end{figure}

\bigskip
We are interested in determining how this similarity differs between all pairs of Pareto-optimal or best known non-dominated solutions (average case), and those solution pairs with the closest objective function values. In our bi-objective setting, the latter subset consists simply of all consecutive pairs of solutions after sorting the Pareto set according to one of the two objectives. We performed these calculations for both sets of instances and for each of the 18 alternative balance objectives. Figure \ref{similarity} presents the overall frequency distributions of common edges between (a) all pairs of solutions in each Pareto set, and (b) all pairs of consecutive solutions (``neighboring'' solutions in the objective space).

From the figure we can see that there are very clear differences in solution similarity between the average case and solutions with the most similar objective function values. Considering all pairs of Pareto-optimal solutions, the percentage of common edges between pairs follows approximately a normal distribution with a mean of around 55\% to 65\% of shared edges, depending on the size of the instances. In contrast, this distribution is highly skewed toward more edges in common when we consider only the subset of ``neighboring'' Pareto-optimal solutions, i.e.~those closest to each other in the objective space. These pairs of solutions have a median similarity of around 75\% to 80\% edges in common.%, likewise greater on the larger instances.

These observations are in line with previous research on fitness landscape analysis for the single-objective CVRP \citep{kubiak2007}, and they reinforce some of the methodological approaches proposed for solving bi-objective VRPs. Clearly there is much potential for exploiting this type of solution similarity within optimization methods. From a managerial perspective, our results suggest that better balanced solutions do not necessarily represent significant departures from the structure of low-cost solutions, and hence are likely to find acceptance in practical settings, especially in light of the generally low marginal cost of balance observed previously.

%% file: sections/5-conclusion.tex
\section{Conclusion}
\label{5-conclusion}

%Workload equity and balanced resource utilization have become increasingly important issues in real-world logistics systems. This stems from the recognition that logistics is not exclusively cost-driven, and that the added cost of more balanced operations can be more than offset by gains through lower overtime hours, improved employee satisfaction, better customer service, and more efficient use of available capacity. The relevance of these considerations is reflected in the number and variety of application cases described in the literature, including service technician routing, vendor-managed inventory systems, waste collection, parcel delivery, and even volunteer organizations. Given their heterogeneous nature, it is not surprising that application-oriented articles have proposed many different forms of ``balance'' objectives. However, the implications of choosing one such objective over another remain largely unexplored. The theory-oriented literature has thus far not adequately addressed this question, having focused mainly on balance in terms of tour length.

In this article, we generalized previous analyses of balanced VRPs to alternative workload resources encountered in practice, and then extended the scope of those analyses to examine further questions of managerial and methodological importance. By considering workload resources which had thus far been unexamined, our study provides more comprehensive insights into the many potential types of balance objectives, and the implications of choosing one such objective over another. Since the VRPs encountered in practice are many times larger than the small instances upon which previous studies have been based, we also replicated our analysis on larger instances to determine the extent to which our conclusions also hold under those circumstances.

Using the analysis of \cite{matl2017a} as a point of departure, we categorized balance objectives according to (a) the workload \textit{resource} which is to be balanced, and (b) the equity \textit{function} used to quantify the degree of balance. Workload resources are either constant-sum or variable-sum, depending on whether the total workload to be distributed is the same for all solutions, or not. Equity functions can be classified in various ways, and we identified two characteristics which \textit{directly interact} with the type of workload resource -- monotonicity, and compatibility with the Pigou-Dalton transfer principle. Based on this analysis we point out some general guidelines for formulating a balance objective. In particular, we emphasize the incompatibility between variable-sum resources and non-monotonic functions, and the greater relevance of the Pigou-Dalton transfer principle when considering functions for constant-sum resources.

Despite the generality of those guidelines, objectives and models satisfying them can still differ significantly in ways which do not lend themselves to a purely analytical treatment. We therefore extended our analysis with a numerical study examining all combinations of three potential workload resources and six alternative equity functions. By considering the bi-objective problem of optimizing both cost and balance, we analyzed the types of compromise solutions which would be found in practice with standard constraint-based or weighted-sum approaches. Our observations can be summarized as follows:

\vspace*{0.25cm}
\begin{itemize}
	\item More complex equity functions provide more potential compromise solutions, as does balancing distance compared to load, and load compared to the number of stops per tour. In practice, if balance considerations are handled with constraints, then a larger set of compromise solutions increases the likelihood that a solution close to a specific constraint value actually exists.
	\item The trade-off between cost and balance depends primarily on the workload resource, and not the equity function. However, the marginal cost of balance was found to be low for all combinations of resource and function, and it was observed to decrease with instance size. This reinforces the practicality of pursuing better balance.
	\item For a given resource, most of the solutions identified by an equity function are not unique to that function and found by at least one other alternative. However, no single equity function finds much more than half of all the potentially ``well-balanced'' solutions for a given resource (according to all considered functions). This implies that there is a large number of well-balanced solutions, of which any function identifies only a limited subset.
%\end{itemize}

% New questions.

%\begin{itemize}
	\item Solutions which are well-balanced in terms of one resource are usually not well-balanced in terms of another. However, for the same resource, solutions optimizing one of the examined equity functions tend to be of high quality also for the others, if not exactly the same. This further underlines the importance of choosing the proper workload resource for a given application -- once a resource has been selected, the equity function has a comparatively limited impact.
	\item Solutions with similar objective values tend to exhibit similar solution structure. The pairs of solutions closest to each other in the objective space were found to have a median of 75\% to 80\% of their edges in common. This lends some credibility to multi-objective solution methods which implicitly assume some form of common structure among high-quality solutions.
\end{itemize}
\vspace*{0.25cm}

In light of the somewhat subjective nature of equity, our study emphasizes the importance of selecting the correct resource to balance. As outlined above, we observed that the subsequent choice of the equity function has relatively less impact on the solutions found -- a conclusion which runs counter to the focus of previous studies which considered only the choice of function. However, we have found that balance -- in all the examined forms -- can be significantly improved at only limited extra cost, allowing the non-monetary benefits of balance to be realized in practice.

%However, we do stress that certain types of functions (non-monotonic) are not compatible with certain types of resources (variable-sum), as they lead to unintended optimization results. Overall, we have shown that balance -- in many different forms -- can be significantly improved at only limited extra cost, allowing the non-monetary benefits of balance to be realized in practice.